\documentclass{amsart}

\usepackage{mathtools, % add-on for amsmath
            amssymb, 
            hyperref, % used for urls and autoref
            graphicx,
}

\usepackage{nicematrix} % for matrices with indexed rows and cols
\NiceMatrixOptions{cell-space-limits = 1.5pt} % spacing for math in tables

\usepackage{caption} % to fix spacing for table captions
\usepackage{tikz}
\usetikzlibrary{arrows, positioning, math, quotes}

\makeatletter
\renewcommand{\sectionautorefname}{\S\@gobble}
\renewcommand{\subsectionautorefname}{\S\@gobble}
\renewcommand{\subsubsectionautorefname}{\S\@gobble}
\makeatother

\usepackage{amsthm}
\usepackage{thmtools}
\theoremstyle{plain}
\newtheorem{thm}{Theorem}[section]
\newtheorem{lem}[thm]{Lemma}
\newtheorem{prop}[thm]{Proposition}
\newtheorem{cor}[thm]{Corollary}

\theoremstyle{definition}
\newtheorem{dfn}[thm]{Definition}
\newtheorem{exam}[thm]{Example}

\theoremstyle{remark}
\newtheorem*{rmk}{Remark}

\newenvironment{pf}[1][\proofname]{\proof}{\endproof}

\usepackage{bbm}
\newcommand{\bbone}{\mathbbm{1}}
\newcommand{\CC}{\mathbb{C}}
\newcommand{\QQ}{\mathbb{Q}}
\newcommand{\RR}{\mathbb{R}}
\newcommand{\ZZ}{\mathbb{Z}}
\newcommand{\cD}{\mathcal{D}}
\newcommand{\cM}{\mathcal{M}}

\newcommand{\paren}[1]{\left( #1 \right)}
\newcommand{\set}[1]{\left\{ #1 \right\}}
\newcommand{\abs}[1]{\left\lvert #1 \right\rvert}

\renewcommand{\phi}{\varphi}

\DeclareMathOperator{\Cay}{Cay}

\usepackage[backend=biber, url=false, doi=false]{biblatex}
\title[]{Graphical Designs find Combinatorial Structures}
\author[]{Zawad Chowdhury, Stefan Steinerberger, Rekha R. Thomas}

\begin{document}

\begin{abstract}
Graphical designs are subsets of vertices of a graph that perfectly average a selected set of eigenvectors of the Graph Laplacian.
We show that in highly-structured graphs, graphical designs can coincide with highly structured and well-known combinatorial objects: orthogonal arrays in hypercube graphs, combinatorial block designs and extremizers of the Erd\H{o}s-Ko-Rado theorem in Johnson graphs, and $t$-wise uniform sets of permutations and symmetric subgroups in normal Cayley graphs on the symmetric group. These connections allow tools from spectral graph theory to bear on these combinatorial objects.
We also show that the central vertex in a Mycielskian is an extremely good design and certain designs of the Mycielskian coincide with designs of the original graph. 
\end{abstract}
 
\maketitle

\tableofcontents

\section{Introduction}
\subsection{Graphical Designs}
Let $G=(V,E)$ be a finite, connected, undirected and unweighted graph on $n$ vertices. Any such graph has an associated \emph{Graph Laplacian} $L=D-A$, where $D \in \mathbb{R}^{n \times n}$ is the diagonal matrix with $D_{ii} = \deg(v_i)$, the degree of the $i$th vertex, and $A \in \mathbb{R}^{n \times n}$ is the adjacency matrix of $G$ with 
\[ A_{ij} = \begin{cases} 1 \qquad &\mbox{if}~(i,j) \in E \\ 0 \qquad &\mbox{otherwise.} \end{cases}\]
The matrix $L \in \mathbb{R}^{n \times n}$ is symmetric and positive semidefinite, with eigenvalues $0 = \lambda_1 < \lambda_2 \leq \lambda_3 \leq \dots \leq \lambda_n$, 
and associated eigenvectors $\phi_1, \dots, \phi_n \in \mathbb{R}^n$. This spectral information captures the underlying geometry of the graph $G$, and a collection of $n$ independent eigenvectors form a `canonical' basis for functions $f:V \rightarrow \mathbb{R}$. These facts form the core of Spectral Graph Theory \cite{chung1997spectral, nica2018intro,grigoryan2018intro}. 

\begin{figure}[!h]
    \begin{minipage}{0.48\textwidth}
        \centering
        \begin{tikzpicture}[scale=0.6]
            \node at (-4, 0) {}; % pushes figure to the right
            \tikzstyle{every node}=[circle, fill=black, inner sep=2pt]
            \foreach \y[count=\a] in {10,9,4}
                {\pgfmathtruncatemacro{\kn}{120*\a-90}
                    \node at (\kn:3) (b\a){} ;}
                \foreach \y[count=\a] in {8,7,2}
                {\pgfmathtruncatemacro{\kn}{120*\a-90}
                    \node at (\kn:1.8) (d\a) {};}
                \foreach \y[count=\a] in {1,5,6}
                {\pgfmathtruncatemacro{\jn}{120*\a-30}
                    \node at (\jn:1.6) (a\a) {};}
                \foreach \y[count=\a] in {3,11,12}
                {\pgfmathtruncatemacro{\jn}{120*\a-30}
                    \node at (\jn:3) (c\a) {};}
            \draw[dashed] (a1)--(a2)--(a3)--(a1);
            \draw[thick] (d1)--(d2)--(d3)--(d1);
            \foreach \a in {1,2,3}
                {\draw[dashed] (a\a)--(c\a);
                \draw[thick] (d\a)--(b\a);}
                \draw[thick] (c1)--(b1)--(c3)--(b3)--(c2)--(b2)--(c1);
                \draw[thick] (c1)--(d1)--(c3)--(d3)--(c2)--(d2)--(c1);
                \draw[dashed] (b1)--(a1)--(b2)--(a2)--(b3)--(a3)--(b1);
        \end{tikzpicture}
    \end{minipage}
        \begin{minipage}{.48\textwidth}
        \centering
        \begin{tikzpicture}[scale=0.6]
            \tikzstyle{every node}=[circle, fill=white, draw=black, inner sep=2pt]

            % graph vertices
            \foreach \a in {1,2,...,9}{
                \node (o\a) at (\a*360/9: 3) {};
            }
            \node (i1) at (0, 1) {};
            \node (i2) at (-0.8, -0.4) {};
            \node (i3) at (0.8, -0.4) {};

            % design vertices
            \foreach \a in {1,5,9}{
                \node[fill=black] at (o\a) {};
            }
            \node[fill=black] at (i2) {};

            % graph edges
            \foreach \a[evaluate={\b=int(\a + 1)}] in {1,2,...,8}{
                \draw (o\a) -- (o\b);
            }
            \draw (o9) -- (o1);
            \foreach \a[evaluate={\c=int(\a + 2)}] in {1,4,7}{
                \draw (o\a) -- (o\c);
            }
            \draw (i1) -- (i2) -- (i3) -- (i1);
            \foreach \a[evaluate={\d=int(3*\a - 1)}] in {1,2,3}{
                \draw (i\a) -- (o\d);
            }
        \end{tikzpicture}
    \end{minipage}
    \vspace{-5pt} 
    \caption{Left: the icosahedron, a spherical 5-design of $\mathbb{S}^2$, is automatically recovered by the spectral geometry. Right: a graphical design with four vertices averaging 11 of 12 eigenvectors in the truncated tetrahedral graph.}
\end{figure}
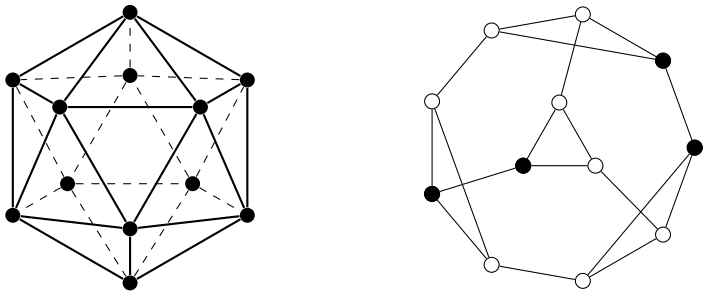

Motivated by \emph{spherical designs} \cite{goethals1977sphericaldesigns,sobolev1962cubature}, one may define \cite{steinerberger2020designs}  a \textbf{graphical design} $\cD \subseteq V$ to be a subset of the vertices of $G$ with the property that, for a certain collection of eigenvectors $\Phi \subseteq \set{\phi_1, \dots, \phi_n}$,
\[
    \frac{1}{|\mathcal{D}|}\sum_{v \in \mathcal{D}} \phi(v) =  \frac{1}{|V|} \sum_{v \in V} \phi(v) \qquad \forall~\phi \in \Phi.
\]
There is a natural question of which subset $\Phi$ one should choose. 
The goal of this paper is to show that, for several particular choices of $\Phi$, well-known combinatorial structures correspond to graphical designs of well-known graphs.
This connection brings a new point of view to the study of these objects. 
 
\subsection{Known Results on Graphical Designs} The main point suggested in \cite{steinerberger2020designs} is that graphical designs, whenever they exist, should be able to capture the underlying geometry of the graph. 
It was also shown in \cite{steinerberger2020designs} that any graphical design that averages a large number of eigenvectors must be well-distributed in the graph, a result motivated by the continuous setting \cite{steinerberger2021spherical}. The first structural result on graphical designs is due to Golubev \cite{golubev2020extremal}: whenever the Hoffman bound or the Cheeger bound is sharp, then the associated set of vertices is a high-quality graphical design. This result has a strong connection to combinatorial structures: it shows that families of subsets attaining the Erd\H{o}s-Ko-Rado bound \cite{erdos1961intersection} are graphical designs of the Kneser graph and maximal intersecting families of permutations in the Deza-Frankl Theorem \cite{deza1977frankl} are graphical designs of the derangement graph. Babecki \cite{babecki2021codes, babecki2022whatis} established a connection to coding theory: the Hamming code is an effective graphical design. Zhu \cite{zhu2023bch} improved the result  using BCH codes.
Babecki and Thomas \cite{babecki2023galeduality} used Gale duality to show that positively weighted
graphical designs of regular graphs are in bijection with the faces of generalized eigenpolytopes of the graph. Steinerberger and Thomas \cite{steinerberger2025equidistribution} showed that graphical designs can be seen as sets with incredibly fast mixing properties when subjected to a random walk. Regarding computation \cite{babecki2024sparsedesigns}, Babecki and Shiroma \cite{babecki2024universality} showed that it is strongly NP-complete to determine if there is a design with size smaller than a fixed upper bound, and NP-hard to find the smallest designs. 

\subsection{Outline of Results}
Graphical designs are based on the eigenfunctions of the graph Laplacian and are thus intrinsically analytic in nature. We show that they also have the ability to recognize structured combinatorial objects within structured graphs. The purpose of this paper is to showcase this combinatorial power of graphical designs in four different settings.
\begin{enumerate}
    \item Graphical designs of the hypercube graph on $\set{0,1}^n$, and more generally in Hamming graphs (\autoref{sec:hamming-section}), are naturally related to \emph{orthogonal arrays}. As a special byproduct, we show that graphical designs for a specific eigenspace correspond to \emph{Hadamard matrices}. Smaller Hamming cubes are graphical designs when the order of eigenspaces is reversed.
    \item Graphical designs of the Johnson graph (\autoref{sec:johnson-section}) correspond to  \emph{combinatorial block designs}. Moreover, reversing the order of eigenspaces, we recover extremal configurations from the Erd\H{o}s-Ko-Rado Theorem.
    \item Moving to the symmetric group $S_n$ (\autoref{sec:S_n-section}), we show that \emph{t-wise uniform sets of permutations} are graphical designs. Reversing the order of the eigenspaces leads to graphical designs that are given by symmetric subgroups.
    \item Finally, we investigate the \emph{Mycielskian} (\autoref{sec:mycielskian-section}) and show that the centrally added vertex plays a special role in Mycielskians of regular graphs: it is a graphical design averaging all but three eigenvectors of the adjacency matrix. The Mycielskian also inherits designs from the original graph.
\end{enumerate}
In \S 2-\S 5 we focus on presenting our results; proofs are in \autoref{sec:proofs}. Our main message is that the graph Laplacian has the ability to find interesting and important combinatorial structures, via the theory of graphical designs of certain structured graphs. 
This is philosophically well aligned with the underlying motivation: graphical designs are, in a certain sense, the optimal samples of the vertices of a graph, while these combinatorial structures (orthogonal arrays, $t$-designs, etc.) naturally arise in statistics as experimental setups that produce a certain optimal sampling of data.

\section{Designs of the Hamming Graph} \label{sec:hamming-section}

We begin by establishing several results about graphical designs of Hamming graphs, using three different orderings of Laplacian eigenvectors. First, in \autoref{sec:laplacian-preliminaries}, we provide background on graphical designs needed throughout the paper. In \autoref{sec:hypercube-designs} we show that the graphical designs in \emph{Laplacian order} are \emph{orthogonal arrays}; the proof involves the Radon transform on a finite group. In the particular case of the second simplest designs of hypercube graphs, we show that the smallest designs correspond to \emph{Hadamard matrices}. In \autoref{sec:reverse-order-designs-hamming} we discuss the \emph{reverse Laplacian order}, in which smaller Hamming cubes are graphical designs. Finally in \autoref{sec:extremal-design-random-walks} we consider the \emph{random walk order} coming from $AD^{-1}$, and resolve a question of Babecki and Thomas about extremal designs of the hypercube graph. All proofs are in \autoref{sec:hamming proofs}.

\subsection{Preliminaries: Graphical Designs and Eigenvector Ordering} \label{sec:laplacian-preliminaries}

\begin{dfn} \label{dfn:graphical-design}
We say that a subset $\cD \subset V$ \textbf{averages} an eigenvector $\phi$ if  
\[\frac{1}{|\mathcal{D}|}\sum_{v \in \mathcal{D}} \phi(v) =  \frac{1}{|V|} \sum_{v \in V} \phi(v).\]
The subsets which average a collection of eigenvectors $\Phi \subset \set{\phi_1, \dots, \phi_n}$ are called \textbf{$\Phi$-designs} of $G$. We also refer to them as graphical designs of $G$ averaging $\Phi$.
\end{dfn}

Let $\bbone_n$ denote the all-ones vector in $\RR^n$. If $G$ is a graph on $n$ vertices, then $\bbone_n$ will be the first eigenvector $\phi_1$ of the graph Laplacian of $G$, with eigenvalue 0. Thus the constant functions $V \to \RR$ comprise the first eigenspace of the Laplacian. 

For any $\phi: V \to \RR$, we have $\sum_{v \in V} \phi(v) = \phi \cdot \bbone_n$.
If $\phi$ is a non-constant eigenvector of the Laplacian, it is orthogonal to $\phi_1 = \bbone_n$, and the sum $\sum_{v \in V} \phi(v)$ is $0$. Therefore a design $\cD$ averages a nonconstant eigenvector $\phi$ if and only if 
\[ \frac1{|\cD|}\sum_{v \in \cD} \phi(v) = 0 \iff \phi\cdot \bbone_{\cD} = 0. \] 
If $\phi$ is the constant eigenvector $\phi_1$, it is averaged by all subsets $\cD \subset V$. Thus we only consider nonconstant eigenvectors when specifying a collection $\Phi$ of eigenvectors to average. We can say $\cD$ averages $\Phi$, or $\cD$ is a $\Phi$-design, if 

\[
\forall~\phi \in \Phi  \qquad  \sum_{v \in \mathcal{D}} \phi(v) =  0. 
\]

We call a $\Phi$-design $\cD$ \emph{minimal} if there is no other $\Phi$-design $\cD'$ with $\cD' \subset \cD$. All $\Phi$-designs are unions of the minimal $\Phi$-designs. We will often be interested in the \emph{smallest} (or minimum) $\Phi$-designs: the designs with the smallest size.

If an eigenbasis of $L$ consists of $\bbone_n = \varphi_1, \varphi_2, \ldots, \varphi_n$, corresponding to eigenvalues $0 = \lambda_1 < \lambda_2 \leq \lambda_3 \leq \cdots \leq \lambda_n$, then two common choices for $\Phi$ are listed below.
\begin{enumerate}
    \item \textbf{Laplacian order} in which $\Phi = \{\phi_2, \ldots, \phi_k\}$ (skipping $\bbone_n = \phi_1$). 
    \item \textbf{Reverse Laplacian order} in which 
    $\Phi = \{ \varphi_n, \varphi_{n-1}, \ldots, \varphi_{k} \}$.
\end{enumerate}

The designs in Laplacian order are the best samples to capture the ``global'' structure of the graph \cite{linderman2020numerint}.
The Laplacian eigenvalue $\lambda_k$ measures how much the corresponding eigenvector $\phi_k$ fluctuates across edges; for normalized $\phi_k$, 
\[ \lambda_k = \sum_{(i,j) \in E} (\phi_k(i) - \phi_k(j))^2.\]
Eigenvectors with small eigenvalues are ``smooth'': they change as little as possible across each edge. Therefore the designs averaging in Laplacian order, starting from the smallest eigenvalues, average the smooth functions $V \to \RR$ and thus philosphically capture the global smooth structure of the graph. On the other hand, designs in reverse Laplacian order average the highly oscillatory eigenvectors with large eigenvalues. Therefore, these designs philosophically capture local structure and are often highly localized subsets of the graph.

We describe graphical designs in Laplacian orders in \S2-3. In \S4, we abandon the Laplacian orders and instead average eigenvectors in an ordering coming from the representation theory of $S_n$. Finally in \S5 we abandon the Laplacian entirely, and instead look at designs averaging eigenvectors of the adjacency matrix. This progression highlights the need to adapt our techniques to the structures at hand.

\subsection{Orthogonal Arrays and Hadamard Matrices} \label{sec:hypercube-designs}

We start this section with the binary case of hypercube graphs, and then move on to general Hamming graphs.

The \textbf{hypercube graph} $H(n, 2)$ is the graph with $2^n$ vertices given by $0/1$ bit strings of length $n$ and an edge between strings that differ in one coordinate. 
In this paper, we think of the vertex set of $H(n,2)$ as the group $(\ZZ/2\ZZ)^n$ to help our proofs. 
The graph induces a metric on the group: the Hamming distance $|x-y|$ between $x, y \in (\ZZ/2\ZZ)^n$ is the number of coordinates they differ in. The Hamming weight $|x|$ of $x$ is the Hamming distance from $x$ to $0$, or equivalently the number of nonzero coordinates in $x$.

Each element $y \in (\ZZ/2\ZZ)^n$ gives us an eigenvector of the $H(n, 2)$ graph Laplacian, with eigenvalue depending on Hamming weight. 
Let $y \cdot x := \sum_{i=1}^n y_ix_i$ be the dot product of $x, y \in (\ZZ/2\ZZ)^n$.
Given $y$, define the function $\chi_y: (\ZZ/2\ZZ)^n \to \RR$ by
\[\chi_y(x) = (-1)^{y \cdot x}.\]
$\chi_y$ is an eigenvector of the graph Laplacian of $H(n,2)$ with eigenvalue $2|y|$. Moreover, the set $\{ \chi_y: y \in (\ZZ/2\ZZ)^n \}$ is a complete set of pairwise orthogonal eigenvectors of the Laplacian.
In fact, the functions $\chi_y$ are the irreducible characters of the group $(\ZZ/2\ZZ)^n$. $H(n, 2)$ is actually the Cayley graph on $(\ZZ/2\ZZ)^n$ with the standard basis vectors $\set{e_1, e_2, \dots, e_n}$ as connection set, and therefore the characters of the group form a basis of eigenvectors of the graph. More details on Cayley graphs and their eigenvectors are provided in \autoref{sec:normal-cayley-eigenvectors}.

\begin{figure}[!h]
    \begin{minipage}{0.48\textwidth}
        \centering
        \begin{tikzpicture}
            \tikzstyle{every node}=[circle, fill=white, draw=black, inner sep=2pt]
            \foreach \x in {0, 1} {\foreach \y in {0, 1} {\foreach \z in {0, 1} {\foreach \w in {0, 1} {
                \node (v\x\y\z\w) at (-0.2*\x + 1*\y + 1.4*\z + 1*\w, 1.4*\x + 1*\y + 0.2*\z - 1*\w) {};
            }}}}
            \foreach \xyzw in {0000, 1110, 0111, 1001} {
                \node[fill=black] at (v\xyzw) {};
            }
            \foreach \x in {0, 1} {\foreach \y in {0, 1} {\foreach \z in {0, 1} {
                \draw (v\x\y\z0) -- (v\x\y\z1);
                \draw (v\x\y0\z) -- (v\x\y1\z);
                \draw (v\x0\y\z) -- (v\x1\y\z);
                \draw (v0\x\y\z) -- (v1\x\y\z);    
            }}}
        \end{tikzpicture}    
    \end{minipage}
    \begin{minipage}{0.48\textwidth}
        \centering
        \begin{tikzpicture}
            \tikzstyle{every node}=[circle, fill=white, draw=black, inner sep=2pt]
            \foreach \x in {0, 1} {\foreach \y in {0, 1} {\foreach \z in {0, 1} {\foreach \w in {0, 1} {
                \node (v\x\y\z\w) at (-0.2*\x + 1*\y + 1.4*\z + 1*\w, 1.4*\x + 1*\y + 0.2*\z - 1*\w) {};
            }}}}
            \foreach \w in {0, 1} {
                \foreach \xyz in {000, 111} {
                    \node[fill=black] at (v\xyz\w) {};
                }
            }
            \foreach \x in {0, 1} {\foreach \y in {0, 1} {\foreach \z in {0, 1} {
                \draw (v\x\y\z0) -- (v\x\y\z1);
                \draw (v\x\y0\z) -- (v\x\y1\z);
                \draw (v\x0\y\z) -- (v\x1\y\z);
                \draw (v0\x\y\z) -- (v1\x\y\z);
            }}}    
        \end{tikzpicture}    
    \end{minipage}
    \caption{The hypercube graph $H(4, 2)$ with some graphical designs. Left: averaging the the first nontrivial eigenspace. Right: averaging all but the middle eigenspace.}
    \label{fig:q4-designs-example}
\end{figure}
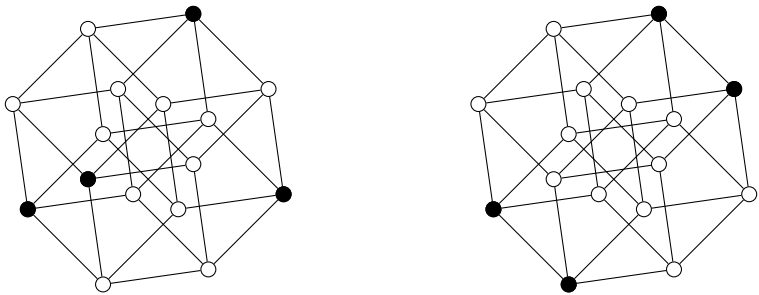

The Laplacian of $H(n,2)$ has $n+1$ eigenspaces, one for each Hamming weight $t = 0, \dots, n$; the eigenspace indexed by $t$ has eigenvalue $2t$ and multiplicity $\binom{n}{t}$. 
To begin with, we assume that designs average entire eigenspaces (see \cite[Lemma 3.2]{babecki2021codes} for motivation). Then our collections $\Phi$ of eigenvectors to average look like
\[\Phi_S := \set{\chi_y: |y| \in S}\]
for some choice of Hamming weights $S \subset [n] = \set{1, 2, \dots, n}$. Hamming weight 0 is ommitted as $\chi_0 = \bbone$ is averaged by all designs; see discussion in \autoref{sec:laplacian-preliminaries}.
For example, 
\[\Phi_{\set{1}} = \set{\chi_y: |y| = 1} = \set{\chi_{e_i}: i \in [n]}.\]
The continuous setting suggests choosing $\Phi$ in \emph{Laplacian order} up to some cutoff, which in our case corresponds to $\Phi_{[t]}$ for some $t \in [n]$.

\begin{exam} \label{exam:cube-design-example}
    Consider the cube graph $H(3, 2)$. The matrix below, whose columns and rows are indexed by $x, y \in (\ZZ/2\ZZ)^3$ respectively, has $\chi_y(x)$ in position $(y,x)$.
    
    \[U = \begin{pNiceMatrix}[first-row, last-col]
        000 & 001 & 010 & 011 & 100 & 101 & 110 & 111 & \\
        1 & 1 & 1 & 1 & 1 & 1 & 1 & 1 & 000\\
        \hline
        1 & -1 & 1 & -1 & 1 & -1 & 1 & -1 & 001\\
        1 & 1 & -1 & -1 & 1 & 1 & -1 & -1 & 010\\
        1 & 1 & 1 & 1 & -1 & -1 & -1 & -1 & 100\\
        \hline 
        1 & -1 & -1 & 1 & 1 & -1 & -1 & 1 & 011\\
        1 & -1 & 1 & -1 & -1 & 1 & -1 & 1 & 101\\
        1 & 1 & -1 & -1 & -1 & -1 & 1 & 1 & 110\\
        \hline
        1 & -1 & -1 & 1 & -1 & 1 & 1 & -1 & 111
    \end{pNiceMatrix}\]
    The rows of this matrix are the eigenvectors $\chi_y$ of $H(3, 2)$. We have grouped them in horizontal strips by eigenvalue/Hamming weight.
    In this case $\Phi_{[1]}$ is the first block of size three, $\set{\chi_{001}, \chi_{010}, \chi_{100}}$. The sets $\cD_1 = \set{000, 111}$ and $\cD_2 = \set{000, 011, 101, 110}$ both average this eigenspace. However, if we average $\Phi_{[2]}$ which includes both blocks of size three, then $\cD_1$ is no longer a design, but $\cD_2$ is. 
\end{exam}
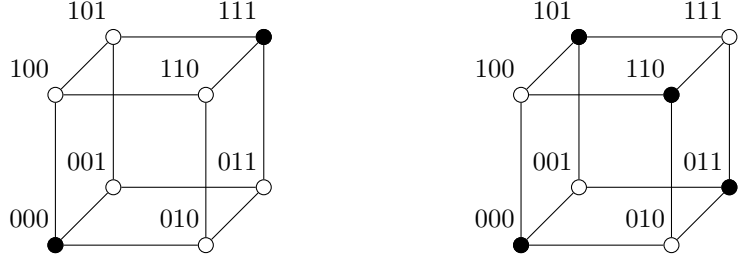
\begin{figure}[ht]
    \begin{minipage}{0.48\textwidth}
        \centering
        \begin{tikzpicture}[scale=2]
            \tikzstyle{every node}=[circle, fill=white, draw=black, inner sep=2pt]
            \foreach \x in {0, 1} {
                \foreach \y in {0, 1} {
                    \foreach \z[evaluate={\l=\y\x\z}] in {0, 1} {
                        \node (v\x\y\z) at (\x, \y, -\z) [label=above left:\l]{};
                    }
                }
            }
            \foreach \xyz in {000, 111} {
                \node[fill=black] at (v\xyz) {};
            }
            \foreach \x in {0, 1} {
                \foreach \y in {0, 1} {
                    \draw (v\x\y0) -- (v\x\y1);
                    \draw (v\x0\y) -- (v\x1\y);
                    \draw (v0\x\y) -- (v1\x\y);
                }
            }
        \end{tikzpicture}
    \end{minipage}
    \begin{minipage}{0.48\textwidth}
        \centering
        \begin{tikzpicture}[scale=2]
            \tikzstyle{every node}=[circle, fill=white, draw=black, inner sep=2pt]
            \foreach \x in {0, 1} {
                \foreach \y in {0, 1} {
                    \foreach \z[evaluate={\l=\y\x\z}] in {0, 1} {
                        \node (v\x\y\z) at (\x, \y, -\z) [label=above left:\l]{};
                    }
                }
            }
            \foreach \xyz in {000, 110, 011, 101} {
                \node[fill=black] at (v\xyz) {};
            }
            \foreach \x in {0, 1} {
                \foreach \y in {0, 1} {
                    \draw (v\x\y0) -- (v\x\y1);
                    \draw (v\x0\y) -- (v\x1\y);
                    \draw (v0\x\y) -- (v1\x\y);
                }
            }
        \end{tikzpicture}
    \end{minipage}
    \caption{Graphical designs of the cube graph $H(3,2)$. Left: the $\Phi_{[1]}$-design $\cD_1$. Right: the $\Phi_{[2]}$-design $\cD_2$.}
    \label{fig:q3-designs-example}
\end{figure}
The bit strings in the design $\cD_2 = \set{000, 011, 101, 110}$ have a curious pattern. Picking any two of the coordinates and deleting the rest, we are left with the set $\set{00, 01, 10, 11}$. In other words, $\cD_2$ has a uniform distribution in all pairs of coordinates. Sets with this property are called \emph{orthogonal arrays}.
\begin{dfn}
    An \textbf{orthogonal array} with parameters $(t, n, q)$ is a set of $q$-ary strings of length $n$ (i.e., a subset of $(\ZZ/q\ZZ)^n$) such that any subset of $t$ coordinates contains every length $t$ string (i.e. every element of $(\ZZ/q\ZZ)^t$) an equal number of times. Such an orthogonal array has \textbf{$q$ levels, $n$ constraints and strength $t$}.
\end{dfn}
In our example above, we see that $\cD_2$ is an orthogonal array with parameters $(2, 3, 2)$. In fact, $\cD_1$ is also an orthogonal array, but with strength 1 instead of 2. This pattern holds in general, and is the main result of this section. We state the result first for hypercube graphs $H(n,2)$; the proof is provided for the general case of Hamming graphs in \autoref{sec:proof-hamming-designs}.

\begin{thm}[Orthogonal Arrays are Designs I] \label{thm:hamming-orthogonalarrays-binary}
    A subset $\cD \subseteq (\ZZ/2\ZZ)^n$ is a $\Phi_{[t]}$-design if and only if $\cD$ is an orthogonal array with parameters $(t, n, 2)$.
\end{thm}

A natural followup question is whether one can construct the smallest $\Phi_{[t]}$-designs for a given $t$. When $t=1$, the smallest designs have size two and are given by antipodal vertices in $H(n, 2)$. When $t=2$, the situation immediately becomes much more complex: the smallest designs correspond to \emph{Hadamard matrices}. 

A matrix in $\{+1, -1\}^{n\times n}$ is called a \textbf{Hadamard matrix} if its rows are pairwise orthogonal. It is well known that an $n \times n$ Hadamard matrix can exist only when $n=1,2$ or is divisible by 4. Proving that such a matrix exists for all $n$ divisible by 4 is the \emph{Hadamard Conjecture}, which has been open for almost a century. The first construction of a Hadamard matrix of order 428 was in 2004 \cite{kharaghani2004hadamard}, and it is still unknown whether there is a Hadamard matrix of order 668 \cite{dokovic2014hadamard}. See \cite[Chapter 2, \S 3]{macwilliams-sloane} for details on the classical results about Hadamard matrices. 

We show that the existence of $\Phi_{[2]}$-designs of $H(n,2)$ with the smallest possible size is equivalent to the Hadamard conjecture. This shows that graphical designs are conceptually hard to find: proving or disproving their existence is, in general, at least as difficult as the Hadamard conjecture, which is very hard indeed.

\begin{thm}[name=Hadamard Matrices are Smallest $\Phi_{[2]}$-Designs, label=thm:hadamard-matrix-designs, restate=thmhadamarddesigns]
    The size of a $\Phi_{[2]}$-design of $H(n, 2)$ is greater than $n$ and divisible by $4$. When $n = 4\ell-1$, a $\Phi_{[2]}$-design with size $4 \ell$ exists if and only if there is a $4\ell \times 4\ell$ Hadamard matrix. 
\end{thm}

Our proof establishes a construction which produces a Hadamard matrix from a $\Phi_{[2]}$-design with size $4\ell$, and vice versa.
For example $\cD_2 = \set{000, 011, 101, 110}$ is a $\Phi_{[2]}$-design of $H(3, 2)$ of size $4\ell$ with $\ell=1$. The corresponding Hadamard matrix produced by the construction outlined in \autoref{sec:proof-hadamard-matrix-design} is
\[\begin{bmatrix}
1 & 1 & 1 & 1 \\
1 & 1 & -1 & -1 \\
1 & -1 & 1 & -1 \\
1 & -1 & -1 & 1
\end{bmatrix}.\]

The story of orthogonal arrays and graphical designs generalizes beautifully if we replace binary with $q$-ary, passing from $H(n,2)$ to the general Hamming graph.
The \textbf{Hamming graph} $H(n, q)$ is the graph with vertex set $(\ZZ/q\ZZ)^n$, with an edge between vertices that differ in one coordinate.
The eigenvectors of the graph Laplacian of $H(n,q)$ are given by the functions $\chi_y: x \mapsto \omega^{y \cdot x}$ for each $y \in (\ZZ/q\ZZ)^n$, where $\omega$ is a primitive $q$-th root of unity. 
 The vectors $\chi_y$ have complex entries if $q \geq 3$, but we can always pick a real basis for the eigenspace which will be a complex linear combination of these vectors.
The notions of Hamming distance and Hamming weight are the same as in $H(n,2)$, and the eigenvector $\chi_y$ has eigenvalue $q|y|$. 
Since the eigenspaces are again indexed by Hamming weight, we consider collections $\Phi_S$, $S \subset [n]$, as before. Once again, $\Phi_{[t]}$-designs in Laplacian order correspond to orthogonal arrays, generalizing \autoref{thm:hamming-orthogonalarrays-binary} to general $q$.

\begin{thm}[name=Orthogonal Arrays are Designs II, label=thm:orth-array-designs, restate=thmortharraydesigns]
    A subset $\cD \subset (\ZZ/q\ZZ)^n$ is a $\Phi_{[t]}$-design if and only if $\cD$ is an orthogonal array with parameters $(t, n, q)$.
\end{thm}

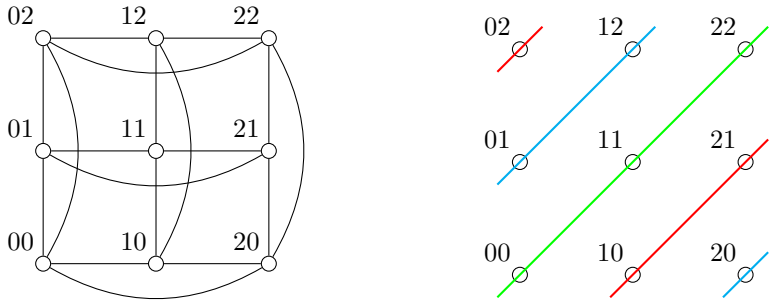
\begin{figure}[!h]
    \centering
    \begin{minipage}{0.48\textwidth}
        \centering
        \begin{tikzpicture}[scale=1.5]
            \tikzstyle{every node}=[circle, fill=white, draw=black, inner sep=2pt]
            \foreach \x in {0, 1, 2} {
                \foreach \y[evaluate={\l=(\x, \y)}] in {0, 1, 2} {
                    \node (v\x\y) at (\x, \y) [label=above left:\l]{};
                }
            }
            \foreach \x in {0, 1, 2} {
                \draw (v\x0) -- (v\x1);
                \draw (v\x1) -- (v\x2);
                \draw (v\x0) to [bend right=30] (v\x2);
            }
            \foreach \y in {0, 1, 2} {
                \draw (v0\y) -- (v1\y);
                \draw (v1\y) -- (v2\y);
                \draw (v0\y) to [bend right=30] (v2\y);
            }

        \end{tikzpicture}
    \end{minipage}
    \begin{minipage}{0.48\textwidth}
        \centering
        \begin{tikzpicture}[scale=1.5]
            \tikzstyle{every node}=[circle, fill=white, draw=black, inner sep=2pt]
            \foreach \x in {0, 1, 2} {
                \foreach \y[evaluate={\l=(\x, \y)}] in {0, 1, 2} {
                    \node (v\x\y) at (\x, \y) [label=above left:\l]{};
                }
            }
            
            \tikzmath{\e = 0.2;}
            \draw[color=red, thick] (0 - \e, 2 - \e) -- (0 + \e, 2 + \e);
            \draw[color=red, thick] (1 - \e, 0 - \e) -- (2 + \e, 1 + \e);
            \draw[color=green, thick] (0 - \e, 0 - \e) -- (2 + \e, 2 + \e);
            \draw[color=cyan, thick] (0 - \e, 1 - \e) -- (1 + \e, 2 + \e);
            \draw[color=cyan, thick] (2 - \e, 0 - \e) -- (2 + \e, 0 + \e);
        \end{tikzpicture}
    \end{minipage}
    \caption{Left: the Hamming Graph $H(2, 3)$. Right: Hyperplanes along one direction in $(\ZZ/3\ZZ)^2$; each color is one hyperplane.}
    \label{fig:hamming-hyperplanes}
\end{figure}

A set $\cD \subset (\ZZ/q\ZZ)^t$ is an orthogonal array whenever certain functions are constant. Whenever $\cD$ is a design, these functions have constant average over all ``hyperplanes'' in $(\ZZ/q\ZZ)^t$ (\autoref{fig:hamming-hyperplanes}). We use the \emph{Radon transform} on finite groups in our proof (\autoref{sec:proof-hamming-designs}) to show that these two formulations are equivalent. 

The fact that $\Phi_{[t]}$-designs are orthogonal arrays lets us deduce two results about the structure of these designs. First, since each of the $q^t$ bit strings in $(\ZZ/q\ZZ)^t$ appear equally many times in $\cD$ for any choice of $t$ coordinates, we get a divisibility constraint on the size of the design.
\begin{cor} \label{cor:divisibility-result}
    If $\cD \subset (\ZZ/q\ZZ)^n$ is a design averaging $\Phi_{[t]}$, then $|\cD|$ is divisible by $q^t$. In particular, we have a bound $|\cD| \ge q^ t$.
\end{cor}
Secondly, an orthogonal array in $(\ZZ/2\ZZ)^n$ still stays an orthogonal array if we get rid of some of its coordinates. Thus we can project designs down. 
\begin{cor}
    Suppose $m \le n$ and $\pi: (\ZZ/q\ZZ)^n \to (\ZZ/q\ZZ)^m$ is a projection map along some choice of $m$ coordinates. If $\cD \subset (\ZZ/q\ZZ)^n$ is a design averaging $\Phi_{[t]}$, then $\pi(\cD) \subset (\ZZ/q\ZZ)^m$ is also a design averaging $\Phi_{[t]}$ in the smaller Hamming graph.
\end{cor}

\subsection{Smaller Hamming Cubes are Reverse Order Designs} \label{sec:reverse-order-designs-hamming}

The Hamming graph has a recursive structure: looking at the subset of vertices with some fixed coordinate values  induces a smaller Hamming graph. These subsets are philosophically opposite to orthogonal arrays; instead of being well distributed amongst all choices of coordinates, they are fixed with respect to a specific choice of coordinates. Quite remarkably, these subsets also serve as graphical designs, averaging eigenvectors in reverse Laplacian order.

\begin{thm}[label=thm:hamming-reverse-designs, restate=thmhammingreversedesigns, name=Hamming Subcubes are Reverse Order Designs]
    Let $\cD$ be a subset of $(\ZZ/q\ZZ)^n$ sharing values in $t$ coordinates, so that the induced subgraph on $\cD$ is isomorphic to the smaller Hamming graph $H(n-t, q)$. Then $\cD$ is a $\Phi_{[n]\setminus [t]}$-design i.e. it averages the first $n-t$ eigenspaces in reverse Laplacian order.
\end{thm}

In the small cases of $H(4, 2)$ and $H(5, 2)$, these designs are the smallest $\Phi_{[n]\setminus[t]}$-designs for all $t$. Computations in those cases suggest that all $\Phi_{[n]\setminus[t]}$-designs have size divisible by $q^{n-t}$, providing a divisibility constraint similar to the one for orthogonal arrays. Proving this divisibility constraint would show that the smaller Hamming cubes are always the smallest designs.

\begin{figure}[!h]
    \begin{minipage}{0.48\textwidth}
        \centering
        \begin{tikzpicture}[scale=2]
            \tikzstyle{every node}=[circle, fill=white, draw=black, inner sep=2pt]
            \foreach \x in {0, 1} {
                \foreach \y in {0, 1} {
                    \foreach \z in {0, 1} {
                        \node (v\x\y\z) at (\x, \y, -\z) {};
                    }
                }
            }
            \foreach \x in {0, 1} {
                \foreach \y in {0, 1} {
                    \node[fill=black] at (v\x\y0) {};
                }
            }
            \foreach \x in {0, 1} {
                \foreach \y in {0, 1} {
                    \draw (v\x\y0) -- (v\x\y1);
                    \draw (v\x0\y) -- (v\x1\y);
                    \draw (v0\x\y) -- (v1\x\y);
                }
            }
        \end{tikzpicture}
    \end{minipage}
    \begin{minipage}{0.48\textwidth}
        \centering
        \begin{tikzpicture}[scale=1]
            \tikzstyle{every node}=[circle, fill=white, draw=black, inner sep=2pt]
            \foreach \x in {0, 1} {\foreach \y in {0, 1} {\foreach \z in {0, 1} {\foreach \w in {0, 1} {
                \node (v\x\y\z\w) at (-0.2*\x + 1*\y + 1.4*\z + 1*\w, 1.4*\x + 1*\y + 0.2*\z - 1*\w) {};
            }}}}
            \foreach \x in {0, 1} {
                \foreach \z in {0, 1} {
                    \node[fill=black] at (v\x0\z0) {};
                }
            }
            \foreach \x in {0, 1} {\foreach \y in {0, 1} {\foreach \z in {0, 1} {
                \draw (v\x\y\z0) -- (v\x\y\z1);
                \draw (v\x\y0\z) -- (v\x\y1\z);
                \draw (v\x0\y\z) -- (v\x1\y\z);
                \draw (v0\x\y\z) -- (v1\x\y\z);    
            }}}
        \end{tikzpicture}    
    \end{minipage}
    \caption{Smaller Hamming cubes are reverse order designs. Left: a $\Phi_{\{2, 3\}}$-design of $H(3, 2)$. Right: a $\Phi_{\{3, 4\}}$-design of $H(4, 2)$.}
    \label{fig:reverse-order-designs}
\end{figure}
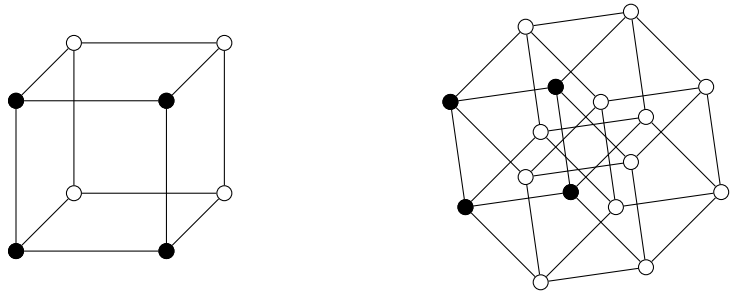

The $t = 1$ case of \autoref{thm:hamming-reverse-designs} was shown by Golubev for $H(n,2)$ \cite[Theorem 3.3]{golubev2020extremal}. The story is actually quite rich in that case: all designs are the union of Hamming subcubes. This fact is not necessarily true for $t > 1$.
\begin{prop}[label=prop:hamming-reverse-eigenpolytope, restate=prophammingreverseeigenpolytope]
    For $i \in [n], a \in \ZZ/q\ZZ$, let $\cD_{i, a}$ be the subset of $(\ZZ/q\ZZ)^n$ consisting of $q$-ary strings $x$ with $x_i = a$. The $qn$ subsets $\cD_{i, a}$ are the minimal designs of $H(n,q)$, averaging $\Phi_{[n]\setminus [1]}$. All other designs are unions of these designs.
\end{prop}

\subsection{Random Walk Order and Extremal Designs} \label{sec:extremal-design-random-walks}

We now turn away from ordering eigenvectors based on the Laplacian and instead look at an ordering motivated by random walks on the graph. For a graph with adjacency matrix $A$ and degree matrix $D$, the matrix $AD^{-1}$ is called the \textbf{random walk matrix} of the graph. This is because for any probability distribution $\pi: V \to \RR$, multiplying the corresponding vector by $AD^{-1}$ gives us the expected probability distribution after one step of the simple symmetric random walk on the graph starting at $\pi$.

If a graph is $d$-regular, $AD^{-1}$ and the Laplacian $L=D-A$ have the same eigenspaces.
There is also a bijection between their eigenvalues. The eigenvalues of $AD^{-1}$ lie in $[-1,1]$ with $1$ corresponding to the first eigenvalue $0$ of $L$. If $L$ has the degree $d$ as an eigenvalue, then $AD^{-1}$ will have $0$ as an eigenvalue. The corresponding eigenvectors are measures on $V$ that decay to zero after one step of the random walk, since multiplying them by $AD^{-1}$ produces the all-zeroes vector.
\begin{figure}[!h]
    \begin{minipage}{0.48\textwidth}
        \centering
        \begin{tikzpicture}[scale=1]
            \tikzstyle{every node}=[circle, fill=white, draw=black, inner sep=2pt]
            \foreach \x in {0, 1} {\foreach \y in {0, 1} {\foreach \z in {0, 1} {\foreach \w in {0, 1} {
                \node (v\x\y\z\w) at (-0.2*\x + 1*\y + 1.4*\z + 1*\w, 1.4*\x + 1*\y + 0.2*\z - 1*\w) {};
            }}}}
            \foreach \zw in {00, 01, 10, 11} {
                \foreach \xy in {00, 11} {
                    \node[fill=red] at (v\xy\zw) {};
                } \foreach \xy in {01, 10} {
                    \node[fill=cyan] at (v\xy\zw) {};
                }
            }
            \foreach \x in {0, 1} {\foreach \y in {0, 1} {\foreach \z in {0, 1} {
                \draw (v\x\y\z0) -- (v\x\y\z1);
                \draw (v\x\y0\z) -- (v\x\y1\z);
                \draw (v\x0\y\z) -- (v\x1\y\z);
                \draw (v0\x\y\z) -- (v1\x\y\z);    
            }}}
        \end{tikzpicture}
    \end{minipage}   
    \begin{minipage}{0.48\textwidth}
        \centering
        \begin{tikzpicture}[scale=1]
            \tikzstyle{every node}=[circle, fill=white, draw=black, inner sep=2pt]
            \foreach \x in {0, 1} {\foreach \y in {0, 1} {\foreach \z in {0, 1} {\foreach \w in {0, 1} {
                \node (v\x\y\z\w) at (-0.2*\x + 1*\y + 1.4*\z + 1*\w, 1.4*\x + 1*\y + 0.2*\z - 1*\w) {};
            }}}}
            \foreach \x in {0, 1} {\foreach \y in {0, 1} {\foreach \z in {0, 1} {
                \draw (v\x\y\z0) -- (v\x\y\z1);
                \draw (v\x\y0\z) -- (v\x\y1\z);
                \draw (v\x0\y\z) -- (v\x1\y\z);
                \draw (v0\x\y\z) -- (v1\x\y\z);
            }}}    
        \end{tikzpicture}    
    \end{minipage}
    \caption{Random walks starting at the eigenvector $\chi_{1100}$ (left), with $AD^{-1}$ eigenvalue $0$, decay to zero after one step (right).}
    \label{fig:random-walk-q4-eigenfunction}
\end{figure}
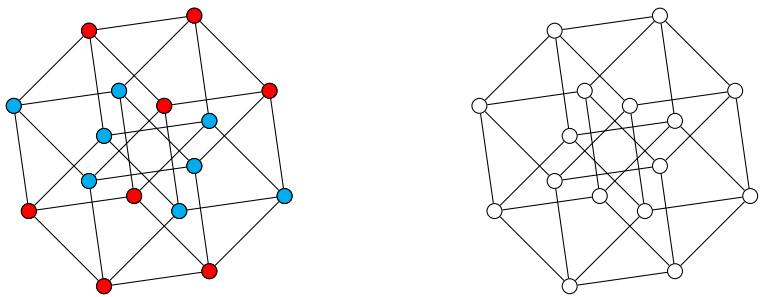

We define the \textbf{random walk order} to be the ordering of eigenvectors in decreasing order of absolute value of $AD^{-1}$ eigenvalue. This corresponds to starting at $1$ at one end of the $AD^{-1}$ spectrum and bouncing back and forth across $0$, advancing towards the middle. This produces an interesting mix of Laplacian and reverse Laplacian orders.

A design in any order is \textbf{extremal} if it averages all but the last eigenspace in that order. Suppose $G = (V, E)$ is a $d$-regular graph with $0$ as an $AD^{-1}$ eigenvalue. Then the extremal designs in random walk order are those which average all eigenspaces besides the one with $AD^{-1}$ eigenvalue 0. We show that these designs are well spread out in the graph.

\begin{prop}[label=prop:random-walk-nbhd, restate=proprandomwalknbhd, name=Extremal Designs Neighbor All Vertices]
    Suppose $G = (V, E)$ is a $d$-regular graph with $0$ as an $AD^{-1}$ eigenvalue, and $\cD \subset V$ an extremal design in random walk order. Then every vertex of $G$ is adjacent to $m$ vertices in $\cD$ for some positive integer $m$. In particular the 1-neighborhood of $\cD$ is the whole graph. 
\end{prop}

\autoref{prop:random-walk-nbhd} explains several of the examples in the original paper on graphical designs by Steinerberger \cite{steinerberger2020designs}, and can be considered a companion to his result about good designs having neighborhoods with exponential volume growth. 
Philosophically, these results further how graphical designs capture the whole graph despite being a small subset of the vertices.

\begin{figure}[!h]
    \begin{minipage}{0.48\textwidth}
        \centering
        \begin{tikzpicture}[scale=1]
            \tikzstyle{every node}=[circle, fill=white, draw=black, inner sep=2pt]
            \foreach \x in {0, 1} {\foreach \y in {0, 1} {\foreach \z in {0, 1} {\foreach \w in {0, 1} {
                \node (v\x\y\z\w) at (-0.2*\x + 1*\y + 1.4*\z + 1*\w, 1.4*\x + 1*\y + 0.2*\z - 1*\w) {};
            }}}}
            \foreach \w in {0, 1} {
                \foreach \xyz in {000, 111} {
                    \node[fill=black] at (v\xyz\w) {};
                }
            }
            \foreach \x in {0, 1} {\foreach \y in {0, 1} {\foreach \z in {0, 1} {
                \draw (v\x\y\z0) -- (v\x\y\z1);
                \draw (v\x\y0\z) -- (v\x\y1\z);
                \draw (v\x0\y\z) -- (v\x1\y\z);
                \draw (v0\x\y\z) -- (v1\x\y\z);    
            }}}
        \end{tikzpicture}
    \end{minipage}    
    \begin{minipage}{0.48\textwidth}
    \centering
        \begin{tikzpicture}[scale=1]
            \tikzstyle{every node}=[circle, fill=black!25, draw=black, inner sep=2pt]
            \foreach \x in {0, 1} {\foreach \y in {0, 1} {\foreach \z in {0, 1} {\foreach \w in {0, 1} {
                \node (v\x\y\z\w) at (-0.2*\x + 1*\y + 1.4*\z + 1*\w, 1.4*\x + 1*\y + 0.2*\z - 1*\w) {};
            }}}}
            \foreach \x in {0, 1} {\foreach \y in {0, 1} {\foreach \z in {0, 1} {
                \draw (v\x\y\z0) -- (v\x\y\z1);
                \draw (v\x\y0\z) -- (v\x\y1\z);
                \draw (v\x0\y\z) -- (v\x1\y\z);
                \draw (v0\x\y\z) -- (v1\x\y\z);
            }}}    
        \end{tikzpicture}    
    \end{minipage}
    \caption{Random walk starting from an extremal design in random walk order (left) immediately equidistributes (right).}
    \label{fig:random-walk-q4-design}
\end{figure}
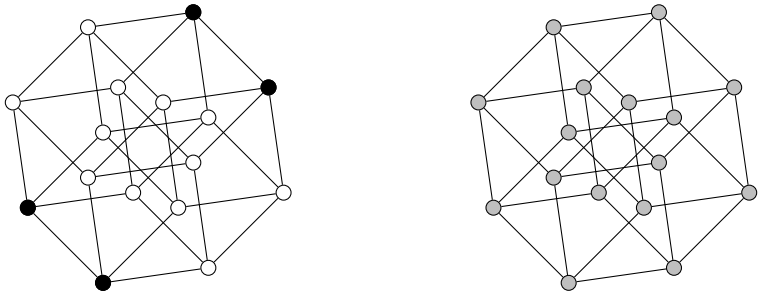

The proof of \autoref{prop:random-walk-nbhd} shows that random walks starting from an extremal design in random walk order immediately equidistributes throughout the graph. This lets us count the adjacencies explicitly and obtain a divisibility bound. In the hypercube graph, this bound answers a question of Babecki and Thomas \cite{babecki2023galeduality}.

\begin{prop}[label=prop:hypercube-extremal-bound, restate=prophypercubeextremalbound]
    For even $n$, if $\cD \subset (\ZZ/2\ZZ)^n$ is an extremal design of $H(n, 2)$ in the random walk order, then $2^n$ divides $n|\cD|$. If $n = 2^tb$ with $b$ odd, the extremal designs must have size at least $2^{n-t}$, and therefore the constructions in \cite[Theorem 5.17]{babecki2023galeduality} are the smallest extremal designs in random walk order.
\end{prop}

\section{Designs of the Johnson graph} \label{sec:johnson-section}

We now consider combinatorial structures in the set of $k$-element subsets of $[n]$, denoted $\binom{[n]}{k}$. If a collection of $k$-element subsets are sufficiently `well distributed' over $[n]$, they form a combinatorial block design (\autoref{dfn:comb-block-design}). These block designs end up being graphical designs of a natural graph on $\binom{[n]}{k}$.

\begin{dfn}
    The \textbf{Johnson graph}, $J(n,k)$ is the graph with vertex set $\binom{[n]}{k}$ and edges between pairs of vertices that share $k-1$ elements.
    Alternatively, if we identify a $k$-element subset of $[n]$ with its indicator vector in $\{0,1\}^n$, then the Johnson graph is the graph whose vertices are the bit strings with Hamming weight $k$, with an edge between vertices at Hamming distance $2$. 
\end{dfn}

The graph $J(n,k)$ is regular of degree $k(n-k)$. Also, $J(n,k)$ is isomorphic to $J(n,n-k)$ and hence we may assume that $k \leq \lfloor n/2 \rfloor$. The octahedral graph $J(4,2)$ is in \autoref{fig:J(4,2)}, while $J(5,2)$ is the complement of the Petersen graph. 
In \autoref{sec:comb-designs} we will see that the graphical designs of $J(n,k)$ in the Laplacian order are \emph{combinatorial block designs}, and in \autoref{sec:reverse-order-designs-johnson} we will see that the designs in reverse Laplacian order are related to the \emph{Erd\"os-Ko-Rado theorem}. 

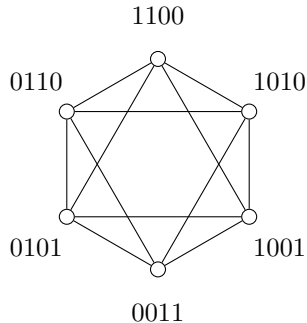
\begin{figure}[!h]
    \centering
    \begin{tikzpicture}[scale=0.7]
        \tikzstyle{every node}=[circle, fill=white, draw=black, inner sep=2pt]

        \node (A) at (90:2) [label=above:1100]{};
        \node (B) at (30:2) [label=above right:1010]{};
        \node (C) at (-30:2) [label=below right:1001]{};
        \node (D) at (-90:2) [label=below:0011]{};
        \node (E) at (-150:2) [label=below left:0101]{};
        \node (F) at (150:2) [label=above left:0110]{};

        \draw (A) -- (B);
        \draw (A) -- (C);
        \draw (A) -- (F);
        \draw (A) -- (E);
        \draw (B) -- (C);
        \draw (F) -- (B);
        \draw (E) -- (F);
        \draw (C) -- (E);
        \draw (D) -- (B);
        \draw (D) -- (C);
        \draw (D) -- (E);
        \draw (D) -- (F);

    \end{tikzpicture}
    \caption{The Johnson Graph $J(4, 2)$}
    \label{fig:J(4,2)}
\end{figure}

\subsection{Combinatorial Block Designs} \label{sec:comb-designs}
The Laplacian of $J(n,k)$ has $k+1$ distinct eigenvalues $\lambda_t$ for $t = 0, \dots, k$ \cite[Theorem 6.3.2]{godsil-meagher}, with multiplicity $\mu_t$ given by
\begin{align} \label{eq:johnson-eigenvals}
    \lambda_t = t (n + 1-t) \qquad \mbox{and} \qquad \mu_t = \frac{n-2t+1}{n-t+1} \binom{n}{t}.
\end{align}
Indeed, the eigenvalues of $J(4,2)$ are $0,4^{(3)},6^{(2)}$ where the exponent denotes multiplicity.
For each eigenvalue  $\lambda_t > 0$, choose an orthonormal basis $\Phi_t$ of its eigenspace.
We are interested in designs averaging $\Phi_S := \cup_{t \in S} \Phi_t$ for $S \subset [k]$. Similar to \autoref{sec:hamming-section}, we consider designs averaging $\Phi_{[t]}$ in \autoref{sec:comb-designs} and designs averaging $\Phi_{[k]\setminus [t]}$ in \autoref{sec:reverse-order-designs-johnson}.

\begin{dfn} \label{dfn:comb-block-design}
    A subset $\cD$ of $\binom{[n]}{k}$ is called a \textbf{$t-(n,k,\lambda)$ design} if every set of $t$ elements in $[n]$ is contained in $\lambda$ elements of $\cD$. The parameters $(n, k, \lambda)$ are often dropped and $\cD$ is referred to as a \textbf{$t$-design}.
\end{dfn}
A $t$-design is also known as a \textbf{combinatorial block design} since the elements of $\cD$ are \emph{blocks} of elements in $[n]$. A $t$-design with $\lambda=1$ is a \textbf{Steiner system}. 
Our first result is that $t$-designs are graphical designs in the Johnson graph.
\begin{thm}[label=thm:comb-designs-are-graphical-designs, restate=thmcombdesigns, name=$t$-designs are Graphical Designs] 
A subset $\cD$ of $\binom{[n]}{k}$ is a $t$-design if and only if it is a $\Phi_{[t]}$-design of the Johnson graph $J(n,k)$.
\end{thm}

The proof is a direct application of a theorem of Delsarte \cite[Theorem 4.7]{delsarte1973thesis} in the context of  \emph{association schemes}. We discuss this connection more broadly in \autoref{sec:association schemes}.

\subsection{Erd\H{o}s-Ko-Rado Constructions are Reverse Order Designs}
\label{sec:reverse-order-designs-johnson} 

Every $t$-element subset of $[n]$ appears equally often in a $t$-design. If instead we pick one $t$-element subset $T$ to be special, and consider all subsets containing $T$, then we get designs in reverse Laplacian order. This is analogous to the story in \autoref{sec:reverse-order-designs-hamming}.

\begin{thm}[label=thm:erdos-ko-rado-designs, restate=thmerdoskoradodesigns, name=Erd\H{o}s-Ko-Rado Extremizers are Reverse Order Designs]
    Given a fixed $t$-element subset $T$ of $[n]$, let $\cD_T \subset \binom{[n]}{k}$ be the set of $k$-element subsets of $[n]$ containing $T$. Then $\cD_T$ averages all $\Phi_\ell$ with $\ell > t$ i.e. $\cD$ is a $\Phi_{[k]\setminus[t]}$-design.
\end{thm}

The generalized Erd\H{o}s-Ko-Rado theorem \cite{wilson1984erdoskorado} states that if $n$ is large enough, the largest size of a collection of $k$-element subsets of $[n]$, in which any two sets intersect in at least $t$ elements, cannot exceed $\binom{n-t}{k-t}$. 
Therefore the designs $\cD_T$ in \autoref{thm:erdos-ko-rado-designs} are the extremal constructions for the generalized Erd\H{o}s-Ko-Rado theorem. That is, the $\cD_T$ are the largest families in $\binom{[n]}{k}$ with pairwise intersections of size at least $t$, when $n$ is sufficiently large.
When $t=1$, these designs are called \textbf{stars}, and they are the largest constructions for the classical Erd\H{o}s-Ko-Rado theorem \cite{erdos1961intersection}. We can show that the stars are not just $\Phi_{[n]\setminus[1]}$-designs, but the minimal designs from which all other designs are constructed.

\begin{prop}[label=prop:johnson-reverse-eigenpolytope, restate=propjohnsonreverseeigenpolytope]
    Let $\cD_i$ be the subset of $\binom{[n]}{k}$ comprised of all $k$-element subsets containing the element $i \in [n]$. The minimal designs of $J(n, k)$ averaging $\Phi_{[k]\setminus [1]}$ are the $\cD_i$ or their complements; all other designs are unions of these designs.
\end{prop}

\subsection{The Connection to Association Schemes} \label{sec:association schemes}

As mentioned already, the proof of \autoref{thm:comb-designs-are-graphical-designs} follows from results on the \emph{Johnson scheme} which is a special case of an \emph{association scheme}. We explain this context briefly, and highlight the key difference between designs of an association scheme and graphical designs. Association schemes are also relevent for some discussion in the next section. See \cite[Chapter 12]{godsil-algcomb}, \cite[Chapter 21]{macwilliams-sloane} or \cite[Chapter 2]{delsarte1973thesis} for more on association schemes.

A \textbf{symmetric association scheme} consists of a set $V$ and $k+1$ graphs $G_0, \ldots, G_k$, each with $V$ as vertex set and disjoint edge sets, such that $G_0$ is the graph of self loops on $V$ and the complete graph on $V$ is the union of the remaining $G_i$. 
For example, the \textbf{Johnson Scheme} is the association scheme on $\binom{[n]}{k}$, identified with vectors of Hamming weight $k$ in $(\ZZ/2\ZZ)^n$, where the edges in the graph $G_i$ are pairs $(x, y)$ with Hamming distance $|x-y| = 2i$. The graph $G_1$ in the Johnson scheme is the Johnson graph $J(n, k)$, and the graph $G_k$ is the \emph{Kneser graph}.

\begin{figure}[!h]
    \begin{minipage}{0.32\textwidth}
        \centering
        \begin{tikzpicture}[scale=0.7]
            \tikzstyle{every node}=[circle, fill=white, draw=black, inner sep=2pt]

            \node (v1) at (90:2) [label=above:1100]{};
            \node (v2) at (30:2) [label=above:1010]{};
            \node (v3) at (-30:2) [label=below:1001]{};
            \node (v4) at (-90:2) [label=below:0011]{};
            \node (v5) at (-150:2) [label=below:0101]{};
            \node (v6) at (150:2) [label=above:0110]{};

            \tikzset{every loop/.style={min distance=15mm}}           
            \foreach \i[evaluate={\a=int(-60*\i)}] in {1,...,6} {
                \draw (v\i) to [in = \a, out = \a-60, loop] ();
            }

        \end{tikzpicture}
    \end{minipage}
    \begin{minipage}{0.32\textwidth}
        \centering
        \begin{tikzpicture}[scale=0.7]
            \tikzstyle{every node}=[circle, fill=white, draw=black, inner sep=2pt]

            \node (A) at (90:2) [label=above:1100]{};
            \node (B) at (30:2) [label=above:1010]{};
            \node (C) at (-30:2) [label=below:1001]{};
            \node (D) at (-90:2) [label=below:0011]{};
            \node (E) at (-150:2) [label=below:0101]{};
            \node (F) at (150:2) [label=above:0110]{};

            \draw (A) -- (B);
            \draw (A) -- (C);
            \draw (A) -- (F);
            \draw (A) -- (E);
            \draw (B) -- (C);
            \draw (F) -- (B);
            \draw (E) -- (F);
            \draw (C) -- (E);
            \draw (D) -- (B);
            \draw (D) -- (C);
            \draw (D) -- (E);
            \draw (D) -- (F);

        \end{tikzpicture}
    \end{minipage}
    \begin{minipage}{0.32\textwidth}
        \centering
        \begin{tikzpicture}[scale=0.7]
            \tikzstyle{every node}=[circle, fill=white, draw=black, inner sep=2pt]

            \node (A) at (90:2) [label=above:1100]{};
            \node (B) at (30:2) [label=above:1010]{};
            \node (C) at (-30:2) [label=below:1001]{};
            \node (D) at (-90:2) [label=below:0011]{};
            \node (E) at (-150:2) [label=below:0101]{};
            \node (F) at (150:2) [label=above:0110]{};

            \draw (A) -- (D);
            \draw (B) -- (E);
            \draw (C) -- (F);
        \end{tikzpicture}
    \end{minipage}
    \caption{The Johnson scheme with $n=4, k=2$. Left: $G_0$, the graph of self loops on $\binom{[4]}{2}$. Middle: $G_1$, the Johnson graph $J(4,2)$. Right: $G_2$, the Kneser graph $K(4,2)$.}
    \label{fig:johnson-scheme}
\end{figure}
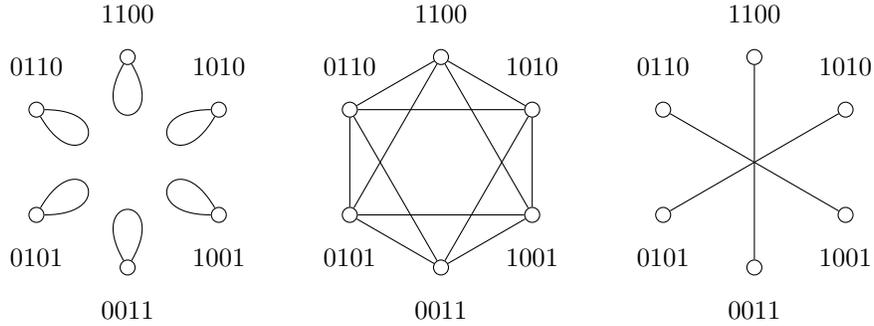

The adjacency matrices $A_i$ of $G_i$ span a vector space in $\RR^{V \times V}$ which happens to be a commutative matrix algebra, called the \emph{Bose-Mesner algebra} of the scheme. This algebra admits another vector space basis of idempotents $\{J_0, \ldots, J_k\}$. Each $J_\ell$ is an orthogonal projector i.e. we can decompose $J_\ell = U_\ell U_\ell^\top$, where $U_\ell$ has pairwise orthogonal columns with norm 1.
Given $\ell \in \set{0, \dots, k}$, the columns of $U_\ell$ will be eigenvectors of the adjacency matrix $A_i$ for all $i$, with all columns of $U_\ell$ having the same eigenvalue. More precisely, the eigenvalues of a fixed $A_i$ are a collection 
of (not necessarily distinct) numbers $p_{i0}, \dots, p_{ik}$ such that  $A_iU_\ell = p_{i\ell}U_\ell$. It is important to note that the numbers $p_{i0}, \dots, p_{ik}$ are not necessarily ordered or distinct. The latter means that an eigenspace of $A_i$ may be the direct sum of column spaces of several $U_\ell$. An extreme example of this ``clumping'' is the following: $A_0 = I$ and hence $1 = p_{00}=p_{01}=\cdots=p_{0k}$ and the eigenspace of $1$ is the direct sum of the columnspaces of all matrices $U_0,\ldots,U_k$. 

For the Johnson scheme, the eigenvalues are given by the \emph{Eberlein polynomials} $E_i(t)$, with $p_{i\ell} = E_i(\ell)$. For example, for $t=0,\dots,k$,
\[E_1(t) = (k-t)(n-k-t) - t\]
are the eigenvalues of the adjacency matrix of the Johnson graph $J(n,k)$. We can verify against \autoref{eq:johnson-eigenvals} that the Laplacian eigenvalues of $J(n,k)$ are then $\lambda_t = k(n-k) - E_1(t)$, where $k(n-k)$ is the degree of $J(n,k)$.

For a symmetric association scheme with $k+1$ graphs and any $T \subset [k]$, Delsarte defined a \textbf{$T$-design} to be a subset $\cD \subset V$ such that $J_\ell\bbone_\cD = 0$ for all $\ell \in T$ \cite[Theorem 3.10]{delsarte1973thesis}. Now, we have
\[J_\ell\bbone_\cD = 0 \iff U_\ell U_\ell^\top \bbone_\cD = 0 \iff \bbone_\cD^\top U_\ell U_\ell^\top \bbone_\cD = 0 \iff U_\ell^\top \bbone_\cD = 0.\]
Thus $\cD$ is a $T$-design of the association scheme if and only if $\cD$ averages all the columns of $U_\ell$ for $\ell \in T$. 

\begin{lem}[label=lem:scheme-to-graphical-design-cond, restate=lemschemetographicaldesigncond]
Suppose there is a graph $G_i$ in the association scheme for which 
\[p_{i0} > p_{i1} > \dots > p_{ik}.\]
Then $\cD$ is a Delsarte $[t]$-design of the scheme if and only if $\cD$ is a $\Phi_{[t]}$-graphical design of $G_i$ in Laplacian order. 
\end{lem}

The proof of \autoref{thm:comb-designs-are-graphical-designs} is precisely to show that the Johnson graph $J(n,k)$ in the Johnson scheme has $k+1$ distinct eigenvalues in the proper order, and thus by \autoref{lem:scheme-to-graphical-design-cond}, $\Phi_{[t]}$-designs of $J(n,k)$ are exactly $[t]$-designs of the Johnson scheme. Then Delsarte's result \cite[Theorem 4.7]{delsarte1973thesis} that $[t]$-designs of the Johnson scheme are combinatorial $t$-designs implies \autoref{thm:comb-designs-are-graphical-designs}. A similar approach with \cite[Theorem 4.4]{delsarte1973thesis} provides a proof of \autoref{thm:orth-array-designs}, different from our proof.

In general, an association scheme is not guaranteed to have any $G_i$ with $k+1$ distinct eigenspaces given by the column spaces of $U_0, \ldots, U_k$. Different sets of column spaces $U_\ell$ may clump together in the different graphs $G_i$. 
This means that there may be no $G_i$ such that for every $t \in [k]$, the $[t]$-designs of the scheme are the $\Phi_{[t]}$-designs of $G_i$. 
A concrete example is the conjugacy class scheme of $S_6$, discussed in \autoref{ex:S6}. 
In such cases, the association scheme designs may average a proper subspace in some eigenspace of every $G_i$, making it impossible to directly compare associaton scheme designs and Laplacian order designs of the graphs in the scheme. 
However, could the Laplacian order designs of the graphs in the scheme have new and interesting  combinatorial structure? In general, association schemes are highly structured combinatorial objects and their designs are also highly structured. Graphical designs, on the other hand, exist for all graphs, and thus provide a more flexible framework for the theory of designs.

\section{Designs of the Symmetric Group} \label{sec:S_n-section}

The combinatorial structures in the symmetric group $S_n$ are revealed by graphical designs of its \emph{normal Cayley graphs}. We describe these graphs and their eigenvectors in \autoref{sec:normal-cayley-eigenvectors}. In \autoref{sec:t-wise-designs} we show that the designs of these graphs correspond to $t$-wise uniform sets of permutations, but following a different order from the Laplacian order. In \autoref{sec:reverse-order-designs-symmetric} we see that reversing this order produces designs which are cosets of symmetric subgroups in $S_n$. 

\subsection{Eigenvectors of Normal Cayley Graphs on $S_n$} \label{sec:normal-cayley-eigenvectors}
There are certain graphs on $S_n$, called normal Cayley graphs, whose eigenvectors and eigenvalues are tied to the representation theory of $S_n$. One example is the \textbf{transposition graph} where the vertex set is $S_n$ and two permutations are connected by an edge if they differ by a transposition. Another example is the \textbf{derangement graph} whose vertex set is $S_n$ and $(\sigma, \pi)$ is an edge if $\sigma(i) \ne \pi(i)$ for all $i \in [n]$.
\begin{dfn}
    Given a group $\Gamma$ and connection set $X \subset \Gamma$, the \textbf{Cayley graph} $\Cay(\Gamma, X)$ is the graph with vertex set $\Gamma$ and edges $(g, gx)$ for $g \in \Gamma, x \in X$. To ensure that the graph is undirected and without loops, $X$ must not contain the identity and be closed under inverses. When $X$ is also closed under conjugation, we say $\Cay(\Gamma, X)$ is \textbf{normal}. 
\end{dfn}

The transposition graph is the Cayley graph on $S_n$ with the transpositions as connection set, and the derangement graph is the Cayley graph whose connection set contains all derangements. Both connection sets are unions of conjugacy classes, and thus both graphs are normal Cayley graphs.

\begin{figure}[!h]
    \begin{minipage}{0.48\textwidth}
        \centering
        \begin{tikzpicture}[scale=1.5]
            \tikzstyle{every node}=[circle, fill=white, draw=black, inner sep=2pt]

            \node (132) at (150:1) [label=above left:(132)]{};
            \node (e) at (90:1) [label=above:e]{};
            \node (123) at (30:1) [label=above right:(123)]{};
            \node (13) at (-30:1) [label=below right:(13)]{};
            \node (12) at (-90:1) [label=below:(12)]{};
            \node (23) at (-150:1) [label=below left:(23)]{};

            \foreach \x in {e, 123, 132} {\foreach \y in {12, 13, 23} {
                \draw (\x) -- (\y);
            }}
        \end{tikzpicture}
    \end{minipage}    
    \begin{minipage}{0.48\textwidth}
        \centering
            \begin{tikzpicture}[scale=1.5]
                \tikzstyle{every node}=[circle, fill=white, draw=black, inner sep=2pt]

                \node (132) at (150:1) [label=above left:(132)]{};
                \node (e) at (90:1) [label=above:e]{};
                \node (123) at (30:1) [label=above right:(123)]{};
                \node (13) at (-30:1) [label=below right:(13)]{};
                \node (12) at (-90:1) [label=below:(12)]{};
                \node (23) at (-150:1) [label=below left:(23)]{};

                \draw (e) -- (123) -- (132) -- (e);
                \draw (12) -- (13) -- (23) -- (12);
            \end{tikzpicture}
    \end{minipage}
    \caption{Normal cayley graphs on $S_3$: transposition graph (left) and derangement graph (right).}
    \label{fig:s3-normal-cayley}
\end{figure}
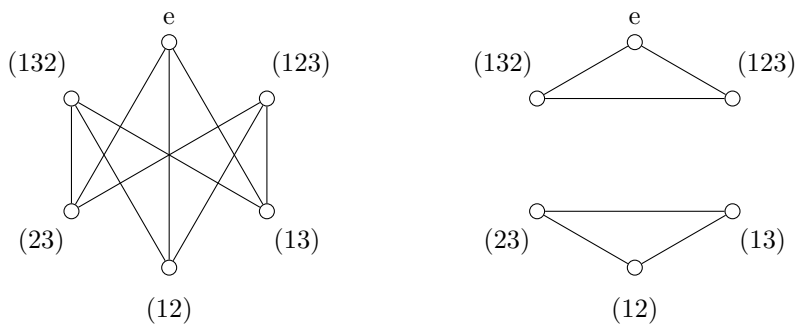

A classical result shows that the entries of irreducible orthogonal representations of a group produce the eigenvectors of its normal Cayley graphs; see \cite[Theorem 10]{ellis2011intersecting} for a proof. In the case of $S_n$, the result shows that:
\begin{itemize}
    \item A partition $p \vdash n$ indexes an irreducible orthogonal representation of $S_n$, $\rho_p: S_n \to \RR^{d_p \times d_p}$, where $d_p$ is the degree of the representation $\rho_p$. 
    \item  Let $\rho_p^{ij}: S_n \to \RR$ be the function $\rho_p^{ij}(\pi) = (\rho_p(\pi))_{i, j}$. If $G = \Cay(S_n, X)$ is a normal Cayley graph, then $\rho_p^{ij}$ is an eigenvector of $G$.
    \item In the graph $\Cay(S_n, X)$, the adjacency matrix eigenvalue of the eigenvector $\rho_p^{ij}$ is $\sum_{x \in X} \chi_p(x) / d_p$, where $\chi_p$ is the character of $\rho_p$.
\end{itemize}

We can check that the eigenvectors $\rho_p^{ij}$ are pairwise orthogonal and there are a total of $n! = |S_n|$ of them. Thus they form an orthogonal basis of eigenvectors for the Laplacian of any normal Cayley graph on $S_n$. Changing the conjugacy classes in the connection set changes the eigenvalues, but the eigenvectors stay the same. 

For any partition $p \vdash n$, let $\Phi_p = \set{\rho_p^{ij}: i,j \in [d_p]}$ be all $d_p^2$ eigenvectors coming from the representation $\rho_p$. For a set $S$ of partitions, define $\Phi_S := \cup_{p \in S} \Phi_p$. 

\begin{exam} \label{exam:S3-eigenvectors}
    Let us construct the eigenvectors for $S_3$, by building \autoref{tab:s3-eigenvectors}. There are three partitions of $3$: $(3), (2, 1)$ and $(1,1,1)$. The corresponding irreducible representations are the trivial, standard and sign representations, respectively. The trivial and sign representations are one dimensional, so we can take their values to form the eigenvectors $\phi_{(3)}^{11}$ and $\phi_{(1,1,1)}^{11}$ respectively.
    The standard representation is two dimensional, with representing matrices
    \begin{alignat*}{2}
    e \mapsto \begin{pmatrix}1 & 0 \\ 0 & 1\end{pmatrix}, \quad &
    (123) \mapsto \begin{pmatrix}-1/2 & -\sqrt{3}/2 \\ \sqrt{3}/2 & -1/2\end{pmatrix},\quad &&
    (132) \mapsto \begin{pmatrix}-1/2 & \sqrt{3}/2 \\ -\sqrt{3}/2 & -1/2\end{pmatrix}, \\
    (23) \mapsto \begin{pmatrix}1 & 0 \\ 0 & -1\end{pmatrix},\quad &
    (12) \mapsto \begin{pmatrix}-1/2 & \sqrt{3}/2 \\ \sqrt{3}/2 & 1/2\end{pmatrix},\quad &&
    (13) \mapsto \begin{pmatrix}-1/2 & -\sqrt{3}/2 \\ -\sqrt{3}/2 & 1/2\end{pmatrix}.
    \end{alignat*}
    Thus the standard representation produces four eigenvectors in $\Phi_{(2, 1)}$, one for each entry of the matrix: $\phi_{(2,1)}^{11}, \phi_{(2,1)}^{12}, \phi_{(2,1)}^{21}, \phi_{(2,1)}^{22}.$

    There are two nontrivial normal Cayley graphs on $S_3$, the transposition graph and the derangement graph (\autoref{fig:s3-normal-cayley}). Let $\lambda_T$ and $\lambda_D$ respectively denote their graph Laplacian eigenvalues for the eigenspaces above. For a given partition $p$, let $\chi_p$ be the corresponding irreducible character. Then the Laplacian eigenvalue for the eigenvectors in $\Phi_p$ can be calculated to be 
    \[\lambda_T = 3 - \frac{3\chi_p((12))}{\chi_p(e)}; \quad \lambda_D = 2 - \frac{2\chi_p((123))}{\chi_p(e)}.\]
    We put all the computed eigenvectors and eigenvalues into a table. 
    \begin{table}[h!]
        $\begin{NiceArray}{|c|c|cccccc|c|c|}
            \hline
            \text{Partition} & \text{Eigenvector} & e & (123) & (132) & (23) & (12) & (13) & \lambda_T & \lambda_D \\
            \hline
            (3) & \rho_{(3)}^{11} & 1 & 1 & 1 & 1 & 1 & 1 & 0 & 0\\
            \hline
            \Block{4-1}{(2,1)} & \rho_{(2,1)}^{11} & 1 & -1/2 & -1/2 & 1 & -1/2 & -1/2 & \Block{4-1}{3} & \Block{4-1}{3} \\
                                & \rho_{(2,1)}^{12} & 0 & -\sqrt{3}/2 & \sqrt{3}/2 & 0 & \sqrt{3}/2 & -\sqrt{3}/2\\
                                & \rho_{(2,1)}^{21} & 0 & \sqrt{3}/2 & -\sqrt{3}/2 & 0 & \sqrt{3}/2 & -\sqrt{3}/2\\
                                & \rho_{(2,1)}^{22} & 1 & -1/2 & -1/2 & -1 & 1/2 & 1/2\\
            \hline
            (1,1,1) & \rho_{(1,1,1)}^{11} & 1 & 1 & 1 & -1 & -1 & -1 & 6 & 0\\
            \hline
        \end{NiceArray}$
        \caption{Eigenvectors of normal Cayley graphs on $S_3$}
        \label{tab:s3-eigenvectors}
    \end{table}

    Note that the derangement graph has the same eigenvalue for $\Phi_{(3)}$ and $\Phi_{(1,1,1)}$. Thus its Laplacian has two eigenspaces rather than three, with $\Phi_{(3)}$ and $\Phi_{(1,1,1)}$ clumped together. This is not uncommon for normal Cayley graphs.
\end{exam}

At this point, it is natural to ask: \emph{which normal Cayley graph on $S_n$ works best for the theory of graphical designs?} The answer, unfortunately, seems to be none in particular, for the usual orderings of eigenvalues. 

\begin{exam} \label{ex:S6} 
In the case of $S_6$, we have 10 different non-identity conjugacy classes. 
None of these conjugacy classes, taken as a connection set for a Cayley graph on $S_6$, avoids the issue of different $\Phi_p$ clumping together in Laplacian eigenspaces. 
We can get around this issue by considering all such graphs at once, using an association scheme. The \emph{conjugacy class scheme of $S_n$} is the symmetric association scheme whose graphs are the normal Cayley graphs with a single conjugacy class as connection set. The $\Phi_p$ are recovered as bases for the eigenspaces of this association scheme. For $S_6$, no graph in the scheme has $11$ distinct eigenvalues, and thus the scheme fails to meet the hypothesis of \autoref{lem:scheme-to-graphical-design-cond}. 
\end{exam} 

In general, the eigenspaces of the conjugacy class scheme of $S_n$, given by the $\Phi_p$, refine the Laplacian eigenspaces of the graphs in the scheme. A proper ordering of the $\Phi_p$ can produce graphical designs which are combinatorial structures, averaging parts of Laplacian eigenspaces rather than whole eigenspaces at once. We discuss two such orders and the corresponding structures in \autoref{sec:t-wise-designs} and \autoref{sec:reverse-order-designs-symmetric}.

\subsection{$t$-wise Uniform Sets and the Fall of Laplacian Order} \label{sec:t-wise-designs}
A well-studied collection of combinatorial structures associated to $S_n$ are $t$-wise uniform sets of permutations; these, like orthogonal arrays and combinatorial block designs, are uniformly distributed among all possibilities at specific coordinates.
\begin{dfn}
    A \textbf{$t$-wise uniform set of permutations} is a subset $\cD \subset S_n$ with a uniform action on any $t$-tuple of elements. In other words, for all choices of distinct $i_1, \dots, i_t \in [n]$ and distinct $j_1, \dots, j_t \in [n]$ we have (for $k=1, \dots, t$),
    \[\frac{\abs{\set{\sigma \in \cD : \sigma(i_k) = j_k}}}{|\cD|} = \frac{\abs{\set{\sigma \in S_n : \sigma(i_k) = j_k}}}{|S_n|} = \frac{(n-t)!}{n!}.\]
    The probability that $\sigma(i_k) = j_k$ for all $1 \le k \le t$ is the same whether $\sigma$ is chosen uniformly from $S_n$ or from $\cD$. Sometimes in short we say $\cD$ is \textbf{$t$-uniform}. 
\end{dfn}
\begin{exam}
    Some examples of $1$-uniform sets of permutations in $S_4$ are
    \[C = \set{e, (1234), (13)(24), (1432)}, \quad V = \set{e, (12)(34), (13)(24), (14)(23)}.\]
    For any $i, j \in [4]$, there is precisely one permutation $\pi \in C$ such that $\pi(i) = j$. The same is true for $V$. Remarkably, they are both subgroups of $S_4$: $C$ is a 4-cycle and $V$ is the Klein four group. Not all $1$-uniform sets are subgroups; we have 
    \[V \cup (12)C = \set{e, (12)(34), (13)(24), (14)(23), (12), (234), (1324), (143)}\]
    a $1$-uniform set in $S_4$ of size 8 that is not a subgroup.

    To be a $2$-wise uniform set in $S_4$, any pair of distinct positions $i_1, i_2 \in [4]$ must be sent to all 12 possible pairs of distinct positions $j_1, j_2 \in [4]$. Thus $2$-uniform sets must have size divisible by 12. If we do not take all of $S_4$, this means our $2$-uniform set will have size exactly 12. Examples of $2$-uniform sets are the alternating group $A_4$ and its complement. One can check that these are all the $2$-uniform sets. 

    For $3$-uniform sets, choosing the images of three numbers in $[4]$ determines the image of the fourth. Thus being $3$-uniform in $S_4$ is the same as being $4$-uniform, and such a set must be $S_4$.
\end{exam}

It is natural to expect, following \autoref{thm:orth-array-designs} and \autoref{thm:comb-designs-are-graphical-designs}, that averaging eigenspaces of some normal Cayley graph on $S_n$ in Laplacian order will produce designs that are $t$-uniform. Unfortunately, this expectation already breaks down for $S_4$: we have three different $t$ for which a set can be $t$-uniform, but four non-constant eigenspaces $\Phi_{(3,1)}$, $\Phi_{(2,2)}$, $\Phi_{(2,1,1)}$, $\Phi_{(1,1,1,1)}$ (recall that $\Phi_{(4)}$ is the eigenspace of constant functions). Looking at which eigenspaces the different $t$-uniform sets actually average, we see that
\begin{center}
    \begin{tabular}{ll}
        $1$-uniform sets average & $\Phi_{(3,1)};$ \\
        $2$-uniform sets average & $\Phi_{(3,1)}, \Phi_{(2,2)}, \Phi_{(2,1,1)};$ \\
        $3$-uniform sets average & $\Phi_{(3,1)}, \Phi_{(2,2)}, \Phi_{(2,1,1)}, \Phi_{(1,1,1,1)}.$
    \end{tabular}
\end{center}
No normal Cayley graph on $S_4$ has eigenspaces that split like this. Looking at the partitions, we spot a different pattern: the new eigenspaces added at each step have the same first part. This pattern works in general.

\begin{thm}[label=thm:t-wise-uniform-designs, restate=twiseuniformdesigns, name=$t$-uniform Sets are Designs]
    A set $\cD \subset S_n$ is $t$-uniform if and only if it is a design averaging all $\Phi_p$ for partitions $p \vdash n$ with first part $p_1 \ge n-t$. 
\end{thm}

\autoref{thm:t-wise-uniform-designs} immediately follows from results by Kuperberg, Lovett and Peled \cite[Section 3.4.2]{kuperberg2017combstructures}; our proof (\autoref{sec:proof-t-uniform-designs}) describes the details to reinterpret their results in the language of graphical designs.

We call the partial order of the partitions implicit in \autoref{thm:t-wise-uniform-designs}, in decreasing order of their first part $p_1$, the \textbf{first part order}. For $S_4$ this order is a coarsening of the transposition graph's Laplacian order. This is not always the case: the transposition graph on $S_6$ has Laplacian eigenvalue 12 for both $\Phi_{(4,1,1)}$ and $\Phi_{(3,3)}$, but the $2$-wise uniform sets on $S_6$ only average the eigenvectors in $\Phi_{(4,1,1)}$ and not necessarily those in $\Phi_{(3,3)}$. Therefore the $t$-wise permutations are designs which do not average entire eigenspaces of the graph Laplacian, but instead distinguish between eigenvectors on the basis of the representation theory of $S_n$. This is quite different from the situation with Hamming and Johnson graphs. We can further show that the first part order never agrees with the Laplacian order of the transposition graph, for sufficiently large $n$. 

\begin{prop}[label=prop:transposition-order, restate=proptranspositionorder]
    For $p \vdash n$, let $\lambda_p$ be the Laplacian eigenvalue of an eigenvector coming from the irreducible representation $\rho_p$ for the transposition graph on $S_n$. 
    If $n \ge 6$, there are two partitions $p, q \vdash n$ with $p_1 > q_1$ ($p < q$ in first part order) but $\lambda_p \ge \lambda_q$ ($p \ge q$ in Laplacian order). Therefore if $n \ge 6$ the first part order does not agree with the Laplacian order.
\end{prop}

The conflict between Laplacian order and first part order is not limited to the transposition graph. \autoref{prop:transposition-order} holds for all normal Cayley graphs on $S_n$ when $n=6,7$. Our computations suggest that the transposition graph is the ``best'' for the result: 
it is the only nontrivial normal Cayley graph on $S_5$ whose Laplacian order agrees with the first part order. However, this conflict poses an interesting question: what are the Laplacian order designs, if not $t$-wise uniform sets? Do they represent a different combinatorial structure? 

\begin{exam}
    In the transposition graph on $S_4$, each $\Phi_p$ has a different eigenvalue. The first two non-trivial eigenspaces in Laplacian order are $\Phi_{(3,1)}$ and $\Phi_{(2,2)}$. A design $\cD \subset S_4$ that averages them is not necessarily a design in first part order, since it may not average $\Phi_{(2,1,1)}$. $\cD$ is 1-uniform as it averages $\Phi_{(3,1)}$, but not necessarily 2-uniform as it does not always average $\Phi_{(2,1,1)}$. The smallest such $\cD$ are of size 12, and up to graph automorphism are either
    \begin{align*}
        &\set{e, (14), (24), (34), (123), (1243), (1423), (1234), (132), (1324), (1342), (1432)},\\
        &\set{e, (1243), (14)(23), (1342), (12), (143), (1324), (234), (13), (243), (1234), (142)}
    \end{align*}
    or $A_4$ (which is also 2-uniform and thus a stronger design). If one can find a special property that distinguishes these $\Phi_{(3,1)} \cup 
    \Phi_{(2,2)}$ designs from other 1-uniform sets in $S_4$, that would be combinatorial information found by the Laplacian order that is not found by the first part order.
\end{exam}

\subsection{Symmetric Subgroups are Reverse Order Designs} \label{sec:reverse-order-designs-symmetric}

The $t$-wise uniform sets of permutations came from averaging eigenvectors in first part order. We might naturally expect that reversing the first part order will yield opposite combinatorial structures, which fix certain positions $\pi(i_k) = j_k$ instead of uniformly distributing every position. If this fixing is done for $n-t$ positions, then the corresponding subsets look like a coset of $S_t$ (under both left and right translation). 

\begin{thm}[label=thm:symm-subgrp-designs, restate=thmsymmsubgrpdesigns, name=Symmetric Subgroups are Reverse Order Designs]
    Suppose $S \subset S_n$ is a subgroup isomorphic to $S_t$ for $t < n$. Then all of its cosets $\cD = \sigma S \sigma'$ are designs of $S_n$ averaging the eigenspaces $\Phi_p$ with $p_1 < n-t$ i.e. in reverse first part order. 
\end{thm}

The generalized \emph{Deza-Frankl conjecture}, proved by Ellis, Friedgut and Pilpel \cite{ellis2011intersecting}, states that if $n$ is sufficiently large compared to $n-t$ and $\cD \subset S_n$ has the property that any two permutations in $\cD$ agree on at least $n-t$ points in $[n]$, then $|\cD| \le t!$. The equality in the bound is achieved precisely by the coset constructions in \autoref{thm:symm-subgrp-designs}. This conjecture generalizes the Erd\H{o}s-Ko-Rado theorem to $S_n$, and thus \autoref{thm:symm-subgrp-designs} is analogous to \autoref{thm:erdos-ko-rado-designs}.

Similar to \autoref{prop:hamming-reverse-eigenpolytope} and \autoref{prop:johnson-reverse-eigenpolytope}, we can fully describe the extremal case of \autoref{thm:symm-subgrp-designs}, when $t=1$ and we are looking at designs averaging all $\Phi_p$ except $\Phi_{(n-1,1)}$. The proofs of those results relied on describing the eigenpolytope of this single omitted eigenspace, and using the Gale duality bijection explained in \autoref{sec:weighted-designs-and-eigenpolytopes}. In the current case, we shall show that the $\Phi_{(n-1, 1)}$-eigenpolytope is the well-known \textbf{Birkhoff polytope}, the convex hull of the permutation matrices in $\RR^{n\times n}$ \cite[Example 0.12]{ziegler1995polytopes}. This is a 0/1-polytope which projects down to the permutahedron, the polytope most commonly associated to $S_n$.

\begin{prop}[label=prop:birkhoff-eigenpolytope, restate=propbirkhoffeigenpolytope]
    The $\Phi_{(n-1, 1)}$-eigenpolytope of normal Cayley graphs on $S_n$ is the Birkhoff polytope. Therefore, the minimal designs averaging the eigenspaces $\Phi_p$ with $p_1 < n-1$, are the subsets of $S_n$ which are cosets of $S_{n-1}$. All other designs averaging these eigenspaces are unions of cosets of $S_{n-1}$.
\end{prop}

\autoref{prop:birkhoff-eigenpolytope} does not hold if $n-1$ is replaced by $n-t$ for $t>1$. The counterexamples provided by Filmus \cite[\S3]{filmus2017comment} are designs averaging all eigenspaces $\Phi_p$ with $p_1 < n-t$, which are not cosets of $S_t$. 

Ellis, Friedgut and Pilpel in their proof of the Deza-Frankl conjecture \cite{ellis2011intersecting} proved results similar to ours, but in the language of Fourier analysis on $S_n$. There are some common elements between our approaches, but the proofs are different.
Their Theorem 7 states that the characteristic functions of the cosets described in \autoref{thm:symm-subgrp-designs} span the space $V_k$ of functions on $S_n$ with Fourier transform supported on irreducible representations corresponding to partitions $p$ with first part $p_1 \ge n-t$. In our language, these are functions orthogonal to eigenfunctions in $\Phi_p$ with $p_1 < n-t$, and thus their theorem can be extended to show the same result as \autoref{thm:symm-subgrp-designs}. 
Their Theorem 28 states that any nonnegative function in the space $V_1$ is a nonnegative combination of the characteristic functions of the minimal designs in \autoref{prop:birkhoff-eigenpolytope}. 
This provides a different proof that the designs in \autoref{prop:birkhoff-eigenpolytope} are the building blocks of the extremal designs.

\section{Designs of the Mycielskian} \label{sec:mycielskian-section}

In this section, we investigate how designs interact with a purely combinatorial construction first given by Jan Mycielski in 1955 \cite{mycielski1955}. Mycielski investigated the notion of whether a graph being triangle-free had any implications for its chromatic number: maybe the chromatic number is determined purely by local obstructions. The \emph{Mycielskian} is a construction that disproves this notion. 
In this section, we show that there is a single vertex in the Mycielskian that averages almost all eigenvectors as a design (\autoref{thm:myc-central-vtx}), and certain designs of the Mycielskian come from designs of the original graph (\autoref{prop:myc-designs-from-original-graph}). 

\begin{figure}[!h]
    \begin{minipage}{0.48\textwidth}
    \centering
    \begin{tikzpicture}[scale=0.8]
        \tikzstyle{every node}=[circle, fill=white, draw=black, inner sep=2pt]
        \def\outr{2.5}       
        \foreach \i in {1,...,5} {
            \node (v\i) at ({18+72*\i}:\outr) [label={18+72*\i}:$v_{\i}$]{};
        }
        \foreach \i[evaluate={\j=int(\i+1)}] in {1,...,4} {
            \draw (v\i) -- (v\j);
        }
        \draw (v5) -- (v1);
    \end{tikzpicture}
    \end{minipage}
    \begin{minipage}{0.48\textwidth}
    \centering
    \begin{tikzpicture}[scale=0.8]
        \tikzstyle{every node}=[circle, fill=white, draw=black, inner sep=2pt]
        \def\outr{2.5}
        \def\inr{1.3}
        
        \foreach \i in {1,...,5} {
            \node (v\i) at ({18+72*\i}:\outr) [label={18+72*\i}:$v_{\i}$]{};
            \node (v'\i) at ({18+72*\i}:\inr) [label={18+72*\i}:$v'_{\i}$]{};
        }
        \node (u) at (0:0) [label=below:$u$]{};
        \foreach \i[evaluate={\j=int(\i+1)}] in {1,...,4} {
            \draw (v'\j) -- (v\i) -- (v\j) -- (v'\i);
        }
        \draw (v'1) -- (v5) -- (v1) -- (v'5);
        \foreach \i in {1,...,5} {
            \draw (v'\i) -- (u);
        }
    \end{tikzpicture}
    \end{minipage}

    \caption{The cycle graph $C_5$ (left) and its Mycielskian, the  Gr\"otzsch graph (right).}
    \label{fig:c5-mycielskian}
\end{figure}
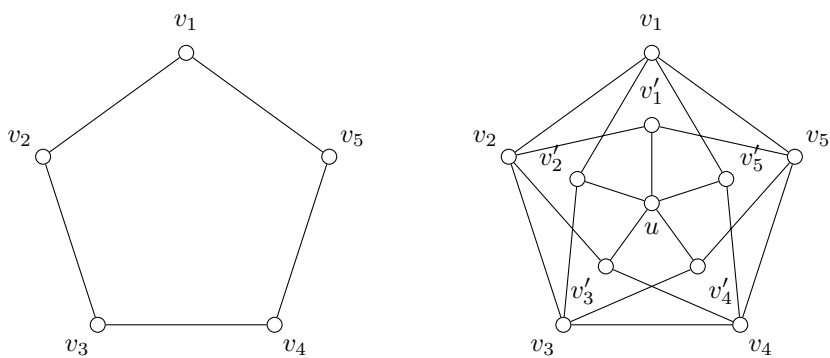

If $G$ is a graph on $n$ vertices $V = \set{v_1, v_2, \dots, v_n}$, the \textbf{Mycielskian $\cM(G)$ of $G$} is a graph on $2n+1$ vertices; the original graph on $v_1, \dots, v_n$ is kept as a subgraph. We add an additional $n$ vertices $V' = \set{v'_1, \dots, v'_n}$, and for each edge $(v_i, v_j) \in E$ we add $(v'_i, v_j)$ and $(v_i, v'_j)$ to the edge set of the Mycielskian. Finally, we add another vertex $u$ and edges connecting it to every $v'_j$ (\autoref{fig:c5-mycielskian}). Mycielski was interested in this construction because starting with a triangle-free graph, the Mycielskian remains triangle-free. However, it is not too difficult to see that the chromatic number of the Mycielskian is larger than the chromatic number of the original graph. Thus iterating the construction produces triangle-free graphs with arbitrarily large chromatic number.

\begin{figure}[!h]
    \centering
    \begin{minipage}{0.48\textwidth}
        \centering
        \includegraphics[width=0.54\textwidth]{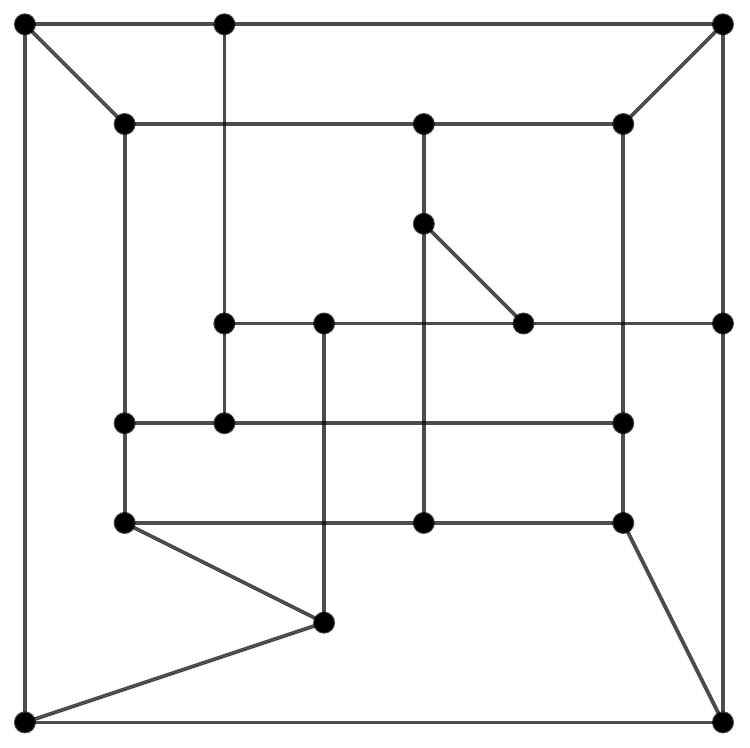}
    \end{minipage}
    \begin{minipage}{0.48\textwidth}
        \centering
        \includegraphics[width=0.66\textwidth]{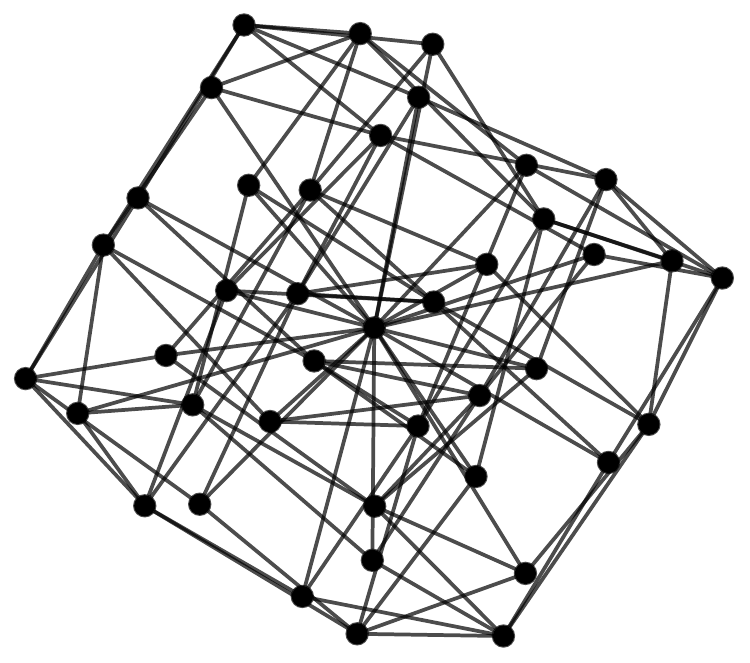}
    \end{minipage}
    \caption{CrossingNumber6B (left) and its Mycielskian (right).}
    \label{fig:crossing-number-6b-mycielskian}
\end{figure}

It is visually clear that the vertex $u$ (connected to all the $v'_i \in V'$) plays a central role in the Mycielskian. This is reflected by the graphical designs. As a concrete example, take $G$ to be the Crossing Number 6B graph on 20 vertices (\autoref{fig:crossing-number-6b-mycielskian}). Then $\cM(G)$ has 41 vertices, and of the 41 eigenvectors of its graph Laplacian all but three (the first, the last and the 18th) have the property that $ \phi(u) = 0.$
Therefore, $\set{u}$ is a small but high-quality graphical design of the Mycielskian, reflecting and emphasizing that $u$ is a central vertex in the Mycielskian construction.

The above example generalizes to any $d$-regular graph $G$ and its Mycielskian. The theory is cleaner if we switch from averaging eigenvectors of the graph Laplacian to averaging eigenvectors of the adjacency matrix. If the adjacency matrix of the original graph $G$ is $A$, then the adjacency matrix of the Mycielskian is
\[
    A(\cM(G)) = \begin{bmatrix}
        A & A & 0\\
        A & 0 & \bbone \\
        0 & \bbone^\top  & 0
    \end{bmatrix}
\] 
in block form. Here $\bbone$ is the all-ones vector in $\RR^n$, where $n$ is the number of vertices in $G$. This block form lets us fully describe the eigenvectors of $A(\cM(G))$, which immediately shows that the central vertex is an extremely good design.

\begin{thm}[label=thm:myc-central-vtx, restate=thmmyccentralvtx, name=Central Vertex is a Good Design]
    Let $G$ be a $d$-regular graph on $n$ vertices. Then $2n-2$ of the $2n+1$ adjacency eigenvectors of $\cM(G)$ vanish at the central vertex $u$, and so $\set{u}$ is a design of $\cM(G)$ averaging all but 3 eigenvectors.
\end{thm}

The block form of $A(\cM(G))$ makes the eigenvalues easy to calculate in terms of the spectrum of $A$ \cite{balakrishnan2012mycielskian}. Suppose the eigenvalues of $A(G)$ are 
\[d = \mu_1 \ge \mu_2 \ge \dots \ge \mu_n.\]
Then $A(\cM(G))$ will have $2n-2$ eigenvalues of the form $\phi\mu_i$ and $\overline{\phi}\mu_i$ for $2 \le i \le n$ (where $\phi, \overline\phi$ are the roots of $x^2 - x - 1$ i.e. the Golden Ratio), and 3 eigenvalues which are the roots of $t^3 - dt^2 - (n+d^2)t + dn = 0$. The $2n-2$ paired eigenvalues $\phi\mu_i$ and $\overline\phi\mu_i$ are the ones with eigenvectors averaged by the central vertex. In fact, designs averaging these eigenvectors come from designs of the original graph.

\begin{prop}[label=prop:myc-designs-from-original-graph, restate=propmycdesignsfromoriginalgraph, name=Mycielskian Designs come from the Original Graph]
    For any $i \in \{2, \dots, n\}$, suppose $\cD \subset V(\cM(G))$ is a design of the Mycielskian averaging the two adjacency eigenvectors with eigenvalues $\phi\mu_i$ and $\overline\phi\mu_i$. Then both $\cD \cap V$ and $\cD \cap V'$ are designs of $G$ averaging the adjacency eigenvector with eigenvalue $\mu_i$.
\end{prop}

In a regular graph, the eigenvalues of the adjacency matrix are in reversed order from the corresponding Laplacian eigenvalues. Now, the Mycielskian $\cM(G)$ is not necessarily regular: the vertices $v_i$, $v'_i$ and $u$ all have different degrees. However, we can use this analogy to motivate the \emph{Mycielskian order} of eigenvectors: going in ascending order of eigenvalues of $A(G)$ ($\mu_n, \dots, \mu_2, \mu_1$ in order), we consider designs averaging the eigenvectors of $A(\cM(G))$ whose eigenvalues come from $\mu_i$ as conjugate pairs (a conjugate triple in the case of $\mu_1$). \autoref{thm:myc-central-vtx} tells us that the central vertex, $\set{u}$, is a design averaging all but 3 eigenvectors. 
That is remarkable: any other designs averaging that many eigenvectors need to cover many vertices.

\begin{prop}[label=prop:myc-design-sizes, restate=propmycdesignsizes]
    Any design averaging the first $2n-2$ eigenvectors in Mycielskian order is either $\set{u}$ or contains at least $n$ vertices. Any design averaging more eigenvectors must average all eigenvectors, and thus contain all $2n+1$ vertices.
\end{prop}
 
The above results can be easily extended further to more general block matrix constructions. Since we are primarily interested in combinatorial structures, we only discuss the Mycielskian with its combinatorial interpretation.
It is interesting that the theory of graphical designs of the Mycielskian shows up so transparently when considering the eigenvectors of the adjacency matrix, but is not as apparent through the graph Laplacian. This suggests that the notion of graphical designs might be of interest for other matrices associated to graphs.

\section{Proofs of Results} \label{sec:proofs}

\subsection{Proofs from \autoref{sec:hamming-section}} \label{sec:hamming proofs}

\subsubsection{Preliminaries on Eigenvectors of $H(n, q)$}
The classical hypercube graph $H(n, 2)$ has a well-known basis of eigenvectors given by $\chi_y(x) = (-1)^{y\cdot x}$ for $y \in (\ZZ/2\ZZ)^n$. It is easy to check that $\chi_y$ is an eigenvector of the Laplacian of $H(n, 2)$ with eigenvalue $2|y|$. There are $2^n$ such eigenvectors and they are pairwise orthogonal, thus they form a full basis of eigenvectors for the Laplacian of $H(n, 2)$.

We can naturally generalize these eigenvectors for the Hamming graph $H(n, q)$, with $\chi_y(x) = \omega^{y \cdot x}$ where $\omega$ is a primitive $q$-th root of unity and $y \in (\ZZ/q\ZZ)^n$. Note that there are $q^n$ such vectors, which is the number of eigenvectors of $H(n, q)$.
If $q \ge 3$, then $\omega$ is a complex number and the vectors $\chi_y$ are complex valued. However, they are still eigenvectors of the graph Laplacian $L$ of $H(n, q)$. 
To check this, we note that the degree of every vertex in $H(n,q)$ is $(q-1)n$ and the vertices adjacent to $x$ are all the elements $z \in (\ZZ/q\ZZ)^n$ at Hamming distance 1 from $x$. Then
\begin{align*}
    (L\chi_y)(x) &= (q-1)n\chi_y(x) - \sum_{|z-x|=1} \chi_y(z) = (q-1)n\omega^{y\cdot x} - \sum_{|z-x|=1} \omega^{y\cdot z} \\
    &= \omega^{y\cdot x}\paren{(q-1)n - \sum_{|z-x|=1} \omega^{y\cdot (z-x)}}.
\end{align*}
The sum over $z-x$ is the same as a sum over all elements with Hamming weight 1 (or, equivalently, Hamming distance 1 from the origin). These take the form $je_i$, where $j \in \set{1, \dots, q-1}$ and $e_i$ is the $i$-th standard basis vector. Now $y\cdot je_i = jy_i$, and so we have
\[\sum_{j=1}^{q-1} \omega^{y \cdot je_i} = \sum_{j=1}^{q-1} \omega^{jy_i} = \begin{cases}
    -1 &\text{ if } y_i \ne 0;\\
    q-1 &\text{ if } y_i = 0.
\end{cases}\]
Let $|y|$ be the Hamming weight of $y$. Then $y_i \ne 0$ for $|y|$ coordinates and $y_i = 0$ for the remaining $n-|y|$. Therefore,
\begin{align*}
    (L\chi_y)(x) &= \omega^{y\cdot x}\paren{(q-1)n - \sum_{i=1}^n \sum_{j=1}^{q-1} \omega^{y\cdot je_i}} \\
    &= \omega^{y\cdot x}\paren{(q-1)n - |y|(-1) - (n-|y|)(q-1)} \\
    &= q|y|\chi_y(x).
\end{align*}

Each $\chi_y$ is thus an eigenvector of the Laplacian of $H(n, q)$ with eigenvalue $q|y|$. There are a total of $q^n$ such eigenvectors and they are pairwise orthogonal, and therefore they form a full complex basis of eigenvectors for the Laplacian of $H(n, q)$. 

The designs of $H(n, q)$ are defined in terms of averaging eigenvectors picked from a real basis. However, our arguments in the coming proofs will simplify if we average the complex basis of eigenvectors $\chi_y$. This does not change the space of designs.

\begin{lem} \label{lem:complex-eigenbasis}
    $\cD \subset (\ZZ/q\ZZ)^n$ is a design averaging all real eigenvectors of the Laplacian of $H(n, q)$ with eigenvalue $qk$ if and only if it averages all the complex eigenvectors $\chi_y$ with $|y| = k$.
\end{lem}
\begin{pf}
    Let $L$ be the graph Laplacian of $H(n, q)$, $N = q^n$ the number of vertices in the graph, and $m = (q-1)^k\binom{n}{k}$ the dimension of the eigenspace of $L$ with eigenvalue $qk$. Suppose $V$ is the complex eigenspace of vectors $v \in \CC^N$ satisfying $Lv = (qk)v$. The vectors $\chi_y$ with $|y| = k$ are a complex basis of $V$.

    Pick a basis $\phi_1, \dots, \phi_m$ of the subspace of real eigenvectors of $L$ with eigenvalue $qk$; we show they are also a complex basis of V. Suppose $\sum w_i \phi_i = 0$ for $w_i \in \CC$. The real and imaginary parts of the $w_i$ are linear dependences of the $\phi_i$ over $\RR$, and thus are all zero. This means the $w_i = 0$ and so the $\phi_i$ are independent over $\CC$. Similarly, for any vector $v \in V$, taking the real and imaginary parts of $v$, we can check that it lies in the complex span of the $\phi_i$. Thus the $\phi_i$ are a basis of $V$.

    A design $\cD$ averages all vectors in a subspace if and only if it averages all vectors in any basis of the subspace. Therefore $\cD$ averaging the basis $\phi_i$ implies that it averages the vectors $\chi_y$, and vice versa.
\end{pf}

\subsubsection{Proof of \autoref{thm:orth-array-designs}} \label{sec:proof-hamming-designs}

How do we check whether $\cD \subset (\ZZ/q\ZZ)^n$ is an orthogonal array? For every choice of $t$ coordinates $I = \{i_1, \dots, i_t\} \subset [n]$, we have a counting function $f_{\cD, I}: (\ZZ/q\ZZ)^t \to \ZZ$, which sends $a \in (\ZZ/q\ZZ)^t$ to the number of elements $x \in \cD$ such that $x_{i_j} = a_j$ for all $j \in [t]$. The set $\cD$ is an orthogonal array with parameters $(t, n, q)$ if and only if $f_{\cD, I}$ is a constant function for all $I$. This is the first step in our proof; the whole argument is structured as follows. 
\begin{figure}[h!]
    \centering
    \begin{tikzpicture}
        \node[rectangle, draw] (n11) at (0, 2) {$\cD$ is an orthogonal array};
        \node[rectangle, draw] (n12) at (6, 2) {$f_{\cD, I}$ is a constant function};
        \node[rectangle, draw] (n21) at (0, 0) {$\cD$ is a $\Phi_{[t]}$-design};
        \node[align=center, rectangle, draw] (n22) at (6, 0) {$f_{\cD, I}$ has constant average \\ over hyperplanes};

        \draw[implies-implies, double equal sign distance] (n11) to[] (n12);
        \draw[implies-implies, double equal sign distance] (n12) to["\autoref{lem:finite-field-radon}"] (n22);
        \draw[implies-implies, double equal sign distance] (n21) to["\autoref{lem:hyperplanes-averaged}"] (n22);
    \end{tikzpicture}
\end{figure}

The other two steps in our proof will involve technical lemmas. Let us first try to establish some motivation for \autoref{lem:hyperplanes-averaged} by proving the equivalence when $q$ is prime, before we state the result in the general case.

A set $\cD$ is a $\Phi_{[t]}$-design if we have $\sum_{x \in \cD} \omega^{y \cdot x} = 0$ for all $y$ with $|y| \in [t]$ (\autoref{lem:complex-eigenbasis}).
Fix one such $y$; suppose $I \subset [n]$ is a set of $t$ coordinates containing the support of $y$.
Let $\pi: (\ZZ/q\ZZ)^n \to (\ZZ/q\ZZ)^t$ be the projection onto the coordinates indexed by $I$. Then we have $y \cdot x = \pi(y) \cdot \pi(x)$, and so to compute $\sum_{x \in \cD} \omega^{y\cdot x} = \sum_{x \in \cD} \omega^{\pi(y) \cdot \pi(x)}$ it is enough to look at the $a \in (\ZZ/q\ZZ)^t$ such that $\pi(x) = a$. 
The value $f_{\cD, I}(a)$ counts the number of such $a$. Hence, for $\Phi_{[t]}$-designs $\cD$,
\[0 = \sum_{x \in \cD} \omega^{y \cdot x} = \sum_{a \in (\ZZ/q\ZZ)^t}\sum_{\substack{x\in \cD,\\ \pi(x) = a}} \omega^{\pi(y) \cdot a} = \sum_{a \in (\ZZ/q\ZZ)^t} f_{\cD, I}(a)\omega^{\pi(y)\cdot a}.\]

Since $\omega$ is a $q$th root of unity, the last sum only cares about the value of $\pi(y) \cdot a$ modulo $q$. Splitting up the sum, we get the equation
\[0 = \sum_{a \in (\ZZ/q\ZZ)^t} f_{\cD, I}(a)\omega^{\pi(y)\cdot a} = \sum_{c=0}^{q-1} \paren{\sum_{\pi(y)\cdot a = c} f_{\cD,I}(a)} \omega^c.\]
Geometrically, we can treat $a \mapsto \pi(y)\cdot a \pmod q$ as a linear function $(\ZZ/q\ZZ)^t \to \ZZ/q\ZZ$, and then its level sets could be considered the \textbf{hyperplanes} of $(\ZZ/q\ZZ)^t$. 
If $q$ is prime, then a weighted sum of $q$-th roots of unity equals zero only if the weights are all equal (this statement is classical \cite{johnsen1985lineare}; it also follows from the argument in \autoref{lem:hyperplanes-averaged}). Thus in this case, we can conclude that the weights, which are the sums of $f_{\cD, I}(a)$ over hyperplanes, are constant.

The same argument works for general $q$ (that is not prime), but we need to take care because some hyperplanes $b \cdot a = c$ would have more points than others.
We call a vector $b \in (\ZZ/q\ZZ)^t$ \textbf{primitive} if $\gcd(b_1, b_2, \dots, b_t, q) = 1.$ Any primitive $b \in (\ZZ/q\ZZ)^t$ and any $c \in \ZZ/q\ZZ$ define a \textbf{primitive hyperplane} $\set{a: b \cdot a = c}$. Every primitive hyperplane will have the correct number $q^{t-1}$ of points. This is because there is some direction $a' \in (\ZZ/q\ZZ)^t$ such that $b\cdot a' \equiv 1 \pmod q$, and by translating in the direction $a'$ we see that the hyperplane $b \cdot a = c$ has the same number of points as the hyperplane $b \cdot a = c+1 = c + b \cdot a'$. Therefore for a fixed primitive $b$, all hyperplanes $b \cdot a = c$ have the same number of points no matter the choice of $c$.

Any nonempty hyperplane will be the disjoint union of primitive hyperplanes. Indeed, if $b$ is not primitive, then pulling out the greatest common divisor of its coordinates and $q$ we can write $b = gb'$ for some primitive $b'$. Now, the left hand side of the equation $b \cdot a = c$ is divisible by $g$. 
If $g$ does not divide $c$, then the equation cannot have any solutions and the hyperplane is empty. Otherwise, we have $c = gc'$, and thus
\[gb'\cdot a = gc' {\pmod q} \iff q \mid g(b'\cdot a - c') \iff \frac{q}{g} \mid b'\cdot a - c'\]
Thus the hyperplane is the disjoint union of the $g$ primitive hyperplanes $\{b' \cdot a = c' + \frac{q}{g}k\}$ for $k = 0, 1, \dots, g-1$, and will have size $gq^{t-1}$. 

\begin{lem}[label=lem:hyperplanes-averaged]
    $\cD \subset (\ZZ/q\ZZ)^n$ is a $\Phi_{[t]}$-design if and only if for all $t$ coordinates $I \subset [n]$, the average of $f_{\cD, I}$ over any nonempty hyperplane in $(\ZZ/q\ZZ)^t$ is $|\cD|/q^t$. 
\end{lem}
\begin{pf}
    Suppose $\cD \subset (\ZZ/q\ZZ)^n$ is a $\Phi_{[t]}$-design. Then $\cD$ averages the eigenvector $\chi_y$ for any $y \ne 0$ with support contained in $t$ coordinates $I \subset [n]$. Therefore,
    \begin{align} \label{eq:prim-hyperplane}
        0 = \sum_{x \in \cD} \omega^{y \cdot x} = \sum_{a \in (\ZZ/q\ZZ)^t} f_{\cD, I}(a)\omega^{\pi(y)\cdot a} = \sum_{c=0}^{q-1} \paren{\sum_{\pi(y)\cdot a = c} f_{\cD,I}(a)} \omega^c.
    \end{align}

    Let us fix $t$ coordinates $I \subset [n]$, and consider hyperplanes coming from $y \in (\ZZ/q\ZZ)^n$ with support contained in $I$. Any hyperplane in $(\ZZ/q\ZZ)^t$ is determined by the equation $\pi(y) \cdot a = c$, for some $y \in (\ZZ/q\ZZ)^n$ with support contained in $I$ and some $c \in \ZZ/q\ZZ$. For this lemma, we need to show that for all such $y, c$ we have
    \[\frac{\sum_{\pi(y)\cdot a = c} f_{\cD, I}(a)}{\abs{\set{a: \pi(y)\cdot a = c}}} = \frac{|\cD|}{q^t}.\]
    It suffices to prove the result only for primitive hyperplanes, as then the average of $f_{\cD, I}$ over any nonempty hyperplane will be the same as the average over the primitive hyperplanes it is a disjoin union of. Primitive hyperplanes always have $q^{t-1}$ points, i.e. $\abs{\set{a: \pi(y)\cdot a = c}} = q^{t-1}$, and so to find the average of $f_{\cD, I}$ it suffices to find its sum over the hyperplane.
    For a fixed primitive $y$, define the function $\sigma: \ZZ/q\ZZ \to \CC$ to be this sum, given by
    \[\sigma(c) = \sum_{\pi(y)\cdot a = c} f_{\cD,I}(a).\]
    Our proof uses the \emph{discrete Fourier transform} $\hat{\sigma}$ of the function $\sigma$. There is a rich theory of such functions, but for our purposes we simply need the following formulas given by the discrete transform and the inverse transform.
    \[\hat{\sigma}(k) = \sum_{c = 0}^{q-1} \sigma(c)(\omega^k)^c; \quad \sigma(c) = \frac{1}{q}\sum_{k=0}^{q-1} \hat{\sigma}(k) (\omega^{-c})^k.\]
    Since $\cD$ is a $\Phi_{[t]}$-design, for any $y$ with support in $I$, \autoref{eq:prim-hyperplane} tells us
    \[0 = \sum_{c=0}^{q-1} \paren{\sum_{\pi(y)\cdot a = c} f_{\cD,I}(a)} \omega^c = \sum_{c=0}^{q-1} \sigma(c)\omega^c = \hat{\sigma}(1).\]
    Now, if $y$ is primitive, then multiplying it by nonzero $k \in (\ZZ/q\ZZ)$ will leave some coordinate nonzero. This multiplication does not change the support, which will still be contained in $I$. Thus \autoref{eq:prim-hyperplane} still applies, and scales to
    \[0 = \sum_{a \in (\ZZ/q\ZZ)^t} f_{\cD, I}(a)\omega^{\pi(ky)\cdot a} = \sum_{c=0}^{q-1} \paren{\sum_{\pi(ky)\cdot a = c} f_{\cD,I}(a)} \omega^{c}. \]
    Let $d = \gcd(k, q)$. Note that $d$ divides the left hand side of the equation $\pi(ky)\cdot a = c$, and so it has solutions modulo $q$ if and only if $d$ also divides $c$. 
    $k$ and $d$ are related by a unit in $(\ZZ/q\ZZ)$, so if $d|c$ we can write $c = kc'$. As $c'$ ranges over $(\ZZ/q\ZZ)$, we will get the same $c=kc'$ appear $d$ different times. Thus we have
    \begin{align} \label{eq:scaled-hyperplanes}
        \sum_{c'=0}^{q-1} \paren{\sum_{\pi(ky)\cdot a = kc'} f_{\cD,I}(a)} \omega^{kc'} = d\sum_{c=0}^{q-1} \paren{\sum_{\pi(ky)\cdot a = c} f_{\cD,I}(a)} \omega^{c} = 0.
    \end{align}
    The first term is $\sum_{c'=0}^{q-1} \sigma(c')\omega^{kc'} = \hat{\sigma}(k)$, and so $\hat{\sigma}$ is zero for all nonzero $k$.
    Finally, $f_{\cD, I}(a)$ counts the number of $x \in \cD$ with $\pi(x) = a$. Summing over all possible $a$, we count all $x \in \cD$ to obtain
    \[\hat{\sigma}(0) = \sum_{c=0}^{q-1} \sigma(c) =  \sum_{a \in (\ZZ/q\ZZ)^t} f_{\cD, I}(a) = |\cD|.\]
    Then, using the inverse transform, we have
    \[\sigma(c) = \frac{1}{q}\sum_{k=0}^{q-1} \hat{\sigma}(k) (\omega^{-c})^k = \frac{\hat{\sigma}(0)}{q} = \frac1{q}|\cD|.\]
    Thus $\sigma(c) = |\cD|/q$ for all $c$, and we conclude that the average of $f_{\cD, I}$ over a primitive hyperplane is the desired $|\cD|/q^t$.

    On the other hand, suppose for any choice of $t$ coordinates $I \subset [n]$ the average of $f_{\cD, I}$ over all nonempty hyperplanes in $(\ZZ/q\ZZ)^t$ is $|\cD|/q^t$. To show that $\cD$ is a $\Phi_{[t]}$-design, it is enough to show that it averages $\chi_y$ for any nonzero $y \in (\ZZ/q\ZZ)^n$ with support contained in some $I$.
    $y$ is not necessarily primitive; pulling out the greatest common divisor of its coordinates and $q$, we have $\pi(y) = gb'$ for some integer $g$ dividing $q$ and some primitive $b' \in (\ZZ/q\ZZ)^t$. Let us consider the function $\sigma$ for hyperplanes with normal vector $b'$, 
    \[\sigma(c) = \sum_{b'\cdot a = c} f_{\cD, I}(a).\]
    The function $\sigma$ is identically $|\cD|/q$ since the average of $f_{\cD, I}$ over any nonempty hyperplane is $|\cD|/q^t$ and the primitive hyperplanes $b' \cdot a = c$ have size $q^{t-1}$ for all $c\in \ZZ/q\ZZ$. Then the discrete Fourier transform $\hat{\sigma}$ is concentrated at 0, and so
    \[0 = \hat{\sigma}(g) = \sum_{c'=0}^{q-1} \sigma(c')\omega^{gc'} = g\sum_{c=0}^{q-1} \paren{\sum_{gb'\cdot a = c}f_{\cD, I}(a)}\omega^c\]
    with the last equation coming from \autoref{eq:scaled-hyperplanes}. Now $gb' = \pi(y)$, so by \autoref{eq:prim-hyperplane} the last term equals $\sum_{x \in \cD} \omega^{y \cdot x}= \sum_{x \in \cD} \chi_y(x)$. Therefore $\cD$ averages $\chi_y$.
\end{pf}

So if $\cD$ is a $\Phi_{[t]}$ design, we know that taking the average of $f_{\cD, I}$ over any hyperplane in $(\ZZ/q\ZZ)^t$ gives us a constant, independent of choice of coordinates $I$ and of hyperplane. If we know the average of a function over all hyperplanes, can we determine the function? This question is answered by the \emph{Radon transform}, and using those ideas lets us prove the following lemma. 
It is the type of statement that one would expect to find in the literature; we have been unable to locate a citation but the lemma may not be a new result.

\begin{lem}[label=lem:finite-field-radon]
    For any function $f: (\ZZ/q\ZZ)^t \to \RR$, the average of $f$ over all nonempty hyperplanes in $(\ZZ/q\ZZ)^t$ is constant if and only if $f$ is a constant function.
\end{lem}

\begin{pf}
    If $f$ is constant, it will clearly have constant average over all hyperplanes. Thus we focus our proof on the other direction. Suppose $C$ is the constant average of $f$ over all nonempty hyperplanes. We shall show that $f$ must be constant.

    Consider the matrix $K: (\ZZ/q\ZZ)^t \times (\ZZ/q\ZZ)^t \to \ZZ$ where $K(a^\star, a)$ is the number of $b \in (\ZZ/q\ZZ)^t$ such that $b \cdot a^* = b \cdot a$. The matrix measures, in a sense, the number of hyperplanes shared by $a^*$ with $a$. 
    If we multiply the matrix $K$ by the vector defined by the function $f$, we get
  \begin{align*}
        (Kf)(a^\star) &= \sum_{a \in (\ZZ/q\ZZ)^t} K(a^\star, a)f(a) \\
        &= \sum_{a \in (\ZZ/q\ZZ)^t} \sum_{\substack{b \in (\ZZ/q\ZZ)^t \\ b\cdot a = b \cdot a^\star}} f(a) 
        = \sum_{b \in (\ZZ/q\ZZ)^t} \sum_{\substack{a \in (\ZZ/q\ZZ)^t \\ b\cdot a = b \cdot a^\star}} f(a).
  \end{align*}
    Summing over $a \in (\ZZ/q\ZZ)^t$ with $b \cdot a = b \cdot a^\star$ is summing over a nonempty hyperplane. Let $g_q(b) := \gcd(b_1, \dots, b_t, q)$; any nonempty hyperplane $b \cdot a = c$ will have $g_q(b)q^{t-1}$ points. Then the sum of $f(a)$ over a nonempty hyperplane will be the average $C$ times the number of points $g_q(b)q^{t-1}$. Thus the above expression simplifies to 
\begin{align*}
    (Kf)(a^\star) &= \sum_{b \in (\ZZ/q\ZZ)^t} \sum_{\substack{a \in (\ZZ/q\ZZ)^t \\ b\cdot a
    = b \cdot a^\star}} f(a) \\
    &= \sum_{b \in (\ZZ/q\ZZ)^t} Cg_q(b)q^{t-1} = Cq^{t-1}\sum_{b \in (\ZZ/q\ZZ)^t} g_q(b).
\end{align*}

    We have shown that $(Kf)(a^\star)$ is a constant independent of $a^\star$, and so $K$ sends the vector $f$ to a multiple of the all-ones vector. Now we will show that the matrix $K$ is invertible, with the all-ones as an eigenvector. This will imply that $K$ sends only constant functions to constant functions, and so $f$ must be constant.

    A vector $b \in (\ZZ/q\ZZ)^t$ is counted in $K(a^\star, a)$ if $b \cdot a^\star = b \cdot a \iff b\cdot (a^\star-a) = 0$. Thus $K(a^\star, a)= g_q(a^\star-a)q^{t-1}$ is counting the points in a nonempty hyperplane. Then $K(a^\star, a)$ is a function of just $a^\star - a$; this implies that the eigenvectors of $K$ will be the characters $\chi_y$ of the group $(\ZZ/q\ZZ)^t$. To check this, we compute
    \begin{align*}
        (K\chi_y)(a^\star) &= \sum_{a \in (\ZZ/q\ZZ)^t} K(a^\star, a)\chi_y(a) = \sum_{a \in (\ZZ/q\ZZ)^t} g_q(a^\star - a)q^{t-1}\omega^{y \cdot a} \\
        &= \omega^{y\cdot a^\star}\sum_{a \in (\ZZ/q\ZZ)^t}g_q(a^\star - a)q^{t-1}\omega^{y\cdot(a - a^\star)} \\
        &= \chi_y(a^\star)\sum_{x \in (\ZZ/q\ZZ)^t} g_q(x)q^{t-1}\omega^{-y\cdot x}.
    \end{align*}
    Therefore $(K\chi_y)(a^\star)$ is $\chi_y(a^\star)$ multiplied by a constant independent of $a^\star$, and so the $\chi_y$ are the eigenvectors of $K$. Since the all-ones vector is $\chi_0$, it is also an eigenvector. To check that $K$ is invertible, we compute the eigenvalues and show that none are zero. We use the identity $g = \sum_{d | g} \phi(d)$. The eigenvalue of $\chi_y$ is
    \begin{align*}
        \sum_{x \in (\ZZ/q\ZZ)^t} g_q(x)q^{t-1}\omega^{-y\cdot x} &= q^{t-1}\sum_{x \in (\ZZ/q\ZZ)^t} \omega^{-y \cdot x} \sum_{d | g_q(x)} \phi(d) \\
        &= q^{t-1}\sum_{d|q} \phi(d) \sum_{\substack{x \in (\ZZ/q\ZZ)^t \\ d | g_q(x)}} \omega^{-y \cdot x}.
    \end{align*}
    Let $q' := q/d$. We have a bijection $x = dx'$ between $x \in (\ZZ/q\ZZ)^t$ with $d | g_q(x)$ and $x' \in (\ZZ/q'\ZZ)^t$. Suppose $\pi$ is the projection map $(\ZZ/q\ZZ)^t \to (\ZZ/q'\ZZ)^t$. Then $\omega^{-y\cdot x} = \omega^{-y \cdot dx'} = \paren{\omega^d}^{- \pi(y) \cdot x'}$; $\omega^d$ is a primitive $q'$-th root of unity. Now, 
    \[\sum_{\substack{x \in (\ZZ/q\ZZ)^t \\ d | g_q(x)}} \omega^{-y \cdot x} = \sum_{x' \in (\ZZ/q'\ZZ)^t} \paren{\omega^d}^{-\pi(y) \cdot x'} = \begin{cases}
        0 & \text{ if } \pi(y) \ne 0 \in (\ZZ/q'\ZZ)^t; \\
        (q')^t & \text{ if } \pi(y) = 0 \in (\ZZ/q'\ZZ)^t.
    \end{cases}\]
    Thus the sum is nonzero only when $q' | g_q(y)$. This simplifies our eigenvalue to 
    \[q^{t-1}\sum_{d|q} \phi(d) \sum_{\substack{x \in (\ZZ/q\ZZ)^t \\ d | g_q(x)}} \omega^{-y \cdot x} = q^{t-1}\sum_{q'|g_q(y)} \phi\paren{q/q'} (q')^t.\]
    $g_q(y)$ is at least 1, so the sum will always have at least one term. All terms in the sum are positive, so we can conclude that the eigenvalue is always positive. Therefore $K$ has nonzero eigenvalues and is thus invertible. $K$ also has the constant vector as an eigenvector, and so $Kf$ being constant implies that $f$ is constant.
\end{pf}

 The term $\sum_{q'|g_q(y)} \phi\paren{q/q'} (q')^t$ is related to Pillai's arithmetic function \cite{pillai1933arithmetic}. Pillai's function is the sum when $t=1$ and $g_q(y) = q$. 
We expect other definitions of $K$ would also be able to prove \autoref{lem:finite-field-radon}; $K(a, a^\star)$ could instead be the number of primitive hyperplanes shared by $a$ and $a^\star$, or a count of shared hyperplanes inversely weighted by the number of points in the hyperplane. These other definitions are more natural, but showing that the resulting matrix is invertible is much harder.

\thmortharraydesigns*
\begin{pf}
    A subset $\cD$ is an orthogonal array with parameters $(t, n, q)$ if and only if $f_{\cD, I}$ is constant for all choices of $t$ coordinates $I \subset [n]$. By \autoref{lem:finite-field-radon}, this is equivalent to $f_{\cD, I}$ having the same average over all hyperplanes in $(\ZZ/q\ZZ)^t$. Finally by \autoref{lem:hyperplanes-averaged}, a subset $\cD$ has constant average over hyperplanes if and only if it is a $\Phi_{[t]}$-design.
\end{pf}

\subsubsection{Proof of \autoref{thm:hadamard-matrix-designs}} \label{sec:proof-hadamard-matrix-design}

We first provide a proof that $\Phi_{[2]}$-designs of $H(n, 2)$ must have size greater than $n$. This lemma is the first part of \autoref{thm:hadamard-matrix-designs}.

\begin{lem}[label=lem:lambda-2-min-design-general-bound]
    If $\cD \subset (\ZZ/2\ZZ)^n$ is a $\Phi_{[2]}$-design, then $|\cD| > n$.
\end{lem}
\begin{pf}
    Since $\cD$ is a $\Phi_{[2]}$-design, by \autoref{cor:divisibility-result} we know $|\cD|$ is divisible by 4. Suppose $|\cD| = 4\ell$. Then for any pair of coordinates $I = (i, j) \subset [n]$, we have $f_{\cD, I}(a) = \ell$ for any possible $a \in (\ZZ/2\ZZ)^2 = \set{00, 01, 10, 11}$. We can apply a graph automorphism switching 1s to 0s, to assume without loss of generality that $0 \in (\ZZ/2\ZZ)^n$ is in $\cD$. Label this as $v_0$, and the other vectors in $\cD$ as $v_1$ through $v_{4\ell - 1}$. Let $z_k$ be the number of zeros in the coordinates of $v_k$; for example $z_0 = n$.

    We now constrain $z_k$ by double counting. Take a count of all pairs of coordinates $(i, j)$ and vectors $v_k$ (with $k \ne 0$) such that the $i$-th coordinate of $v_k$ is 0. If we pick $(i, j)$ first, there are $n(n-1)$ choices of pairs of coordinates, and then $\ell - 1$ vectors $v_k$ where the coordinates are $00$ (removing $v_0$ from the count) and $\ell$ vectors $v_k$ where the coordinates are $01$. On the other hand, if we pick $v_k$ first there are $z_k$ choices of $i$, and then $n-1$ remaining choices for $j$. This gives an equation
    \[(n-1)\sum_{k=1}^{4\ell - 1} z_k = n(n-1)(\ell - 1 + \ell) \quad \implies \quad \sum_{k=1}^{4\ell - 1} z_k = n(2\ell - 1).\]
    
    Now count pairs of coordinates $(i, j)$ and vectors $v_k$ (with $k \ne 0$) such that both the $i$-th and $j$-th coordinates of $v_k$ are zero. There are now $\ell -1$ choices for $v_k$ once you fix $(i, j)$. On the other hand, once you fix $v_k$, you have $z_k(z_k - 1)$ choices for $(i, j)$. This gives us

    \[\sum_{k=1}^{4\ell - 1} z_k(z_k - 1) = n(n-1)(\ell - 1) \]
    which implies that 
    \[\sum_{k=1}^{4\ell - 1} z_k^2 = \sum_{k=1}^{4\ell - 1} z_k(z_k -1) + \sum_{k=1}^{4\ell - 1} z_k = n(n-1)(\ell - 1) + n(2\ell - 1).\]
    Applying the Cauchy-Schwarz inequality, we have
    \[\paren{\sum_{k=1}^{4\ell - 1} z_k^2}\paren{\sum_{k=1}^{4\ell - 1} 1} \ge \paren{\sum_{k=1}^{4\ell - 1} z_k}^2 \]
    which implies  
    \[\paren{n(n-1)(\ell - 1) + n(2\ell - 1)}\paren{4\ell - 1} \ge \paren{n(2\ell-1)}^2.\]
    After some calculation, this simplifies to $4 \ell - 1 \ge n$ and hence $|\cD| = 4 \ell > n$.
\end{pf}
\begin{rmk}
    A different approach to this argument can be used to bound the size of $\Phi_{[t]}$-designs of $H(n, 2)$ for all values of $t$, not just $t=2$. The results needed, along with a construction showing that this bound is optimal up to constants, are described by Alon and Spencer in \cite[Section 16.2]{alon2008probabilistic}.
\end{rmk}

\thmhadamarddesigns*
\begin{pf}
    The connection between Hadamard matrices and strength 2 orthogonal arrays is well known. For example, it is stated as a problem in \cite[Chapter 11.8]{macwilliams-sloane}. We write out the details in this proof.

    \autoref{cor:divisibility-result} shows that the design has size divisible by 4, and \autoref{lem:lambda-2-min-design-general-bound} shows that it has size greater than $n$. Thus we focus on constructing a Hadamard matrix from a design with size $4\ell$ and vice versa.
    Let us start with a size $4\ell$ design $\cD \subset (\ZZ/2\ZZ)^n$ with $n = 4\ell - 1$. Then the Cauchy-Schwarz inequality in the proof of \autoref{lem:lambda-2-min-design-general-bound} is tight, which implies that the $z_k$ are all equal to some $z$. Therefore
    \[n(2\ell - 1) = \sum_{k=1}^{4\ell - 1} z_k = nz \quad \implies \quad z = 2\ell - 1.\]
    Thus the Hamming distance
    \[d(v_0, v_k) = n - z_k = 2\ell.\]
    However, $v_0$ is not special! We can apply a graph automorphism to move any vector $v_k$ to $v_0$ while preserving Hamming distance. This means that any pairwise distance in $\cD$ is $2\ell$. 
    Now we can apply the transformation $0 \mapsto 1, 1 \mapsto -1$ on the coordinates of the vectors in $\cD$, and then add a final $4\ell$-th coordinate of all 1s. This gives us $4\ell$ rows in $\{1, -1\}^{4\ell}$, each pairwise different at $2\ell$ coordinates and thus orthogonal. This is the desired Hadamard matrix.

    If we start with a Hadamard matrix, we can negate columns of the matrix, while keeping rows mutually orthogonal, until one row is all 1s. Then every other row must have an equal number of 1s and -1s to be orthogonal to the first row. Remove the all 1s row, apply the transformation $1 \mapsto 0, -1 \mapsto 1$ coordinatewise, and transpose to get a $4\ell \times (4\ell - 1)$ matrix. 
    Each column of this matrix has an equal number of 1s and 0s. Since the columns came from orthogonal rows, the number of $00$ and $11$ in a pair of columns (which contributed $+1$ to the dot product as Hadamard matrix rows) equals the number of $01$ and $10$ in the pair of columns (which contributed $-1$). Interpreting each row as a vertex in $H(4\ell - 1, 2)$, these two facts establish the set as a strength 2 orthogonal array and thus a $\Phi_{[2]}$-design. This establishes the other direction of the correspondence.
\end{pf}

\subsubsection{Proof of \autoref{thm:hamming-reverse-designs}}

\thmhammingreversedesigns*
\begin{pf}
    These smaller Hamming cubes are produced by fixing a choice of $t$ coordintes $I \subset [n]$ and $a \in (\ZZ/q\ZZ)^t$. Let $\cD = \cD_{I, a}$ be the copy of $H(n-t, q) \subset H(n, q)$ containing the vectors $x \in (\ZZ/q\ZZ)^n$ with $\pi(x) = a$. Let $\bbone_\cD$ be the indicator vector of the subset $\cD$; we shall show that $\bbone_\cD$ is a linear combination of eigenvectors in $\Phi_{[t]}$ along with the constant vector. Since eigenvectors are orthogonal, this means that $\bbone_\cD$ will average all eigenvectors in the remaining eigenspaces $\Phi_{[n]\setminus [t]}$.

    For $b \in (\ZZ/q\ZZ)^t$, we can obtain a lift $\overline{b}\in (\ZZ/q\ZZ)^n$ supported on $I$, taking on the same values as $b$ in $I$ and 0 outside of $I$. In particular, the lift satisfies $\overline{b} \cdot x = b \cdot \pi(x)$. Consider $f: (\ZZ/q\ZZ)^n \to \CC$ defined by
    \[f(x) := \frac1{q^t}\sum_{b \in (\ZZ/q\ZZ)^t} \omega^{-b \cdot a} \chi_{\overline{b}}(x) = \frac1{q^t}\sum_{b \in (\ZZ/q\ZZ)^t} \omega^{b \cdot (\pi(x) - a)}.\]

    If $\pi(x) = a$, the summands are all 1 and thus we get $f(x) = 1$. Otherwise $\pi(x) - a \ne 0$, and so summing over all powers of $\omega$ as $b$ ranges over $(\ZZ/q\ZZ)^t$ yields 0. Therefore $f$ is an indicator for whether $\pi(x) = a$; so if $\cD = \cD_{I, a}$ then $f = \bbone_\cD$. We have thus proved
    \[\bbone_\cD = f = \frac1{q^t}\sum_{b \in (\ZZ/q\ZZ)^t} \omega^{-b \cdot a} \chi_{\overline{b}}\]
    lies in the span of $\set{\bbone, \Phi_{[t]}}$, and so $\cD$ averages $\Phi_{[n]\setminus[t]}$.
\end{pf}

\subsubsection{Weighted Graphical Designs and Eigenpolytopes} \label{sec:weighted-designs-and-eigenpolytopes}

We briefly describe \emph{weighted graphical designs} and their bijection with faces of a \emph{generalized eigenpolytope} of $G$. These are the tools used for proving  \autoref{prop:hamming-reverse-eigenpolytope}, as well as \autoref{prop:johnson-reverse-eigenpolytope} and \autoref{prop:birkhoff-eigenpolytope}.

Suppose $G=(V,E)$ is a regular graph and $\Phi$ a collection of eigenvectors 
of the Laplacian of $G$. A collection of vertices $\cD \subset V$ is a \textbf{positively weighted $\Phi$-design} if there exist positive numbers (weights) $w_v$ such that 
\[    
    \forall \varphi \in \Phi, \qquad \sum_{v \in \cD} w_v \varphi(v) = \frac{1}{|V|} \sum_{v \in V} \varphi(v) = 0.
\]
The $\Phi$-designs considered in this paper (\autoref{dfn:graphical-design}) are \emph{uniformly weighted}, namely they have $w_v = 1$ for all $v \in \cD$. 
It was shown in \cite{babecki2023galeduality}  that positively weighted designs are in bijection with the faces of an eigenpolytope of $G$, via the tool of \emph{Gale duality}. We describe this result, since it is the tool for proving a few of our results starting with \autoref{prop:hamming-reverse-eigenpolytope}.

Let $\bar{\Phi}$ denote the eigenvectors of the Laplacian that complete $\Phi$ to form a full eigenbasis, and let ${U}_{\bar \Phi}$ be the matrix whose rows consist of the eigenvectors in $\bar{\Phi}$. 
The $\bar{\Phi}${\bf -eigenpolytope} of $G$ is the convex hull of the columns of ${U}_{\bar{\Phi}}$. 

\begin{thm}[{\cite[Theorem 3.8]{babecki2023galeduality}}] \label{thm:gale-duality-bijection}  
A subset $\cD$ of $V$ is a positively weighted $\Phi$-design of $G$ if and only if the the columns of ${U}_{\bar{\Phi}}$ indexed by $V \setminus \cD$ are precisely the columns that lie on a face of the $\bar{\Phi}$-eigenpolytope of $G$.
\end{thm}

In particular, minimal positively weighted $\Phi$-designs correspond to facets of the $\bar{\Phi}$-eigenpolytope. 
In extraordinary situations, such as in the proof of 
\autoref{prop:hamming-reverse-eigenpolytope} and other similar propositions throughout this paper, all facets of the relevant eigenpolytope correspond to uniformly weighted designs.

\subsubsection{Proof of \autoref{prop:hamming-reverse-eigenpolytope}}
\prophammingreverseeigenpolytope*
\begin{pf}
    By \autoref{thm:gale-duality-bijection}, the  $\Phi_{[n]\setminus [1]}$ designs of $H(n, q)$ are in bijection  with the faces of the $\Phi_{[1]}$-eigenpolytope of $H(n,q)$. This eigenpolytope is well known \cite[Section 5]{godsil1998eigenpolytope}; its skeleton is the graph $H(n, q)$, and its facets correspond to the $qn$ sets $\set{x \in (\ZZ/q\ZZ)^n: x_i \ne a}$ for $i \in [n], a \in \ZZ/q\ZZ$. The sets $\cD_{i,a}$ are the complements of the facets, and thus the minimal designs averaging $\Phi_{[n]\setminus[1]}$. By \autoref{thm:hamming-reverse-designs} they are uniformly weighted. Any other design is the union of these minimal designs as any face of the eigenpolytope is the intersection of some facets. 
\end{pf}

\subsubsection{Proof of \autoref{prop:random-walk-nbhd}}
\proprandomwalknbhd*
\begin{pf}
    If $\cD \subset V$ is an extremal design in random walk order, we have $\bbone_\cD  = c\bbone + f$ with $\bbone$ the all-ones vector on $V$, $c$ some constant and $f$ an eigenvector with $AD^{-1}$ eigenvalue 0. The other eigenspaces do not contribute to $\bbone_\cD$ as $\cD$ averages them.

    Since $\bbone$ and $f$ are eigenvectors of $AD^{-1}$ with eigenvalue 1 and 0 respectively, the random walks starting at measure $\bbone$ stay constant, while those starting at measure $f$ immediately decay to the 0 measure everywhere. Now
    \begin{align} \label{eq:random-walk-converges}
        AD^{-1}\bbone_\cD = AD^{-1}(c\bbone) + AD^{-1}f = c\bbone.
    \end{align} 
    With the appropriate normalization, this shows that the random walk starting from the uniform measure on $\cD$ converges after a single step to the uniform measure on the whole graph. 
    This means every vertex of $G$ is adjacent to an equal number of vertices $m$ in $\cD$. As $m\ge 1$, the 1-neighborhood of $\cD$ is the whole graph.
\end{pf}

\subsubsection{Proof of \autoref{prop:hypercube-extremal-bound}}

\prophypercubeextremalbound*
\begin{pf}
    If $n$ is even, then $G = H(n, 2)$ is an $n$-regular graph with $n$ as a Laplacian eigenvalue, and thus 0 as an $AD^{-1}$ eigenvalue. If $\cD$ is an extremal design in the random walk order, applying \autoref{prop:random-walk-nbhd} tells us that every vertex in $H(n, 2)$ is adjacent to $m$ vertices in $\cD$. We know that every vertex in $\cD$ is adjacent to exactly $n$ vertices in $H(n, 2)$; from this double count of the adjacencies, we have 
    \[n|\cD| = m|(\ZZ/2\ZZ)^n| = m2^n.\]
    Thus $2^n$ divides $n|\cD|$. The rest of the proposition's conclusions follow. 
\end{pf}

The ideas in \autoref{prop:random-walk-nbhd} and \autoref{prop:hypercube-extremal-bound} can be applied to non-extremal designs. Suppose $\bbone_\cD = c\bbone + f$, where $f$ is a linear combination of eigenvectors with $AD^{-1}$ eigenvalue less than a bound $\lambda$. Then \autoref{eq:random-walk-converges} becomes
\[AD^{-1}\bbone_\cD < c\bbone + \lambda f.\]
This does not give us an exact description of adjacencies and thus cannot produce a divisibility constraint. However, it can provide a lower bound on the size of $\cD$. For example, for the hypercube graph $H(n, 2)$ with $n$ even, we can deduce that if $\cD$ averages all but the last $k$ eigenspaces in random walk order, then 
\[|\cD| > \frac{2^n}{n}\paren{1 - \frac{k^2}{3n}}.\]
In the extremal case of $k = 1$, this almost gives us the divisibility bound $|\cD| \ge 2^n/n$. This approach can still produce bounds for far weaker designs in random walk order, as long as the design averages all eigenspaces with large $AD^{-1}$ eigenvalue.

\subsection{Proofs from \autoref{sec:johnson-section}}

\subsubsection{Proof of \autoref{thm:comb-designs-are-graphical-designs}}
\thmcombdesigns*
\begin{pf}
    By \autoref{lem:scheme-to-graphical-design-cond} it suffices to prove that $p_{10} > p_{11} > \cdots > p_{1k}$. Since the $p_{1i}$ are the eigenvalues of the adjacency matrix of $J(n,k)$, we can instead show that $\lambda_0 < \lambda_1 < \cdots < \lambda_k$ where $\lambda_t = t(n+1-t)$ are the Laplacian eigenvalues of $J(n, k)$. Indeed, $\lambda_t < \lambda_{t+1}$ for $t < k \leq n/2$.
\end{pf}

\subsubsection{Proof of \autoref{thm:erdos-ko-rado-designs}}
\thmerdoskoradodesigns*
\begin{pf}
    Let $B_t$ be the matrix with rows indexed by $\binom{[n]}{k}$ and columns indexed by $\binom{[n]}{t}$, with $B_t(S, T) = 1$ if $T \subset S$ for $S \in \binom{[n]}{k}, T \in \binom{[n]}{t}$, and $B_t(S, T) = 0$ if $T \not\subset S$. The columns of $B_t$ are the indicator vectors of our designs $\cD_T$.

    Now consider the matrix $C_t = B_tB_t^\top$. Delsarte showed in his thesis that the $C_t$ form a basis of the Bose-Mesner algebra of the Johnson scheme, and $C_t$ is in the span of $J_0, \dots, J_t$ (see \cite[Equation 4.29]{delsarte1973thesis}; note that Delsarte's notation for $B_t$ is $A_i$). This means that $J_\ell C_t = 0$ for all $\ell > t$. 
    From the theory of association schemes, $J_\ell = U_\ell U_\ell^\top$ and the columns of $U_\ell$ generate the eigenspace $\Phi_\ell$. Then
    \[J_\ell C_t = 0 \implies U_\ell U_\ell^\top B_t B_t^\top  = 0 \implies U_\ell^\top B_t = 0.\]
    
    The columns of $B_t$ are $\bbone_\cD$ for a design $\cD = \cD_T$, as $T$ ranges over $\binom{[n]}{t}$. Therefore $\bbone_\cD$ is perpendicular to all the eigenvectors in $\Phi_\ell$ for $\ell > t$, and we can conclude that $\cD$ averages $\Phi_{[k]\setminus [t]}$.
\end{pf}

\subsubsection{Proof of \autoref{prop:johnson-reverse-eigenpolytope}}
\propjohnsonreverseeigenpolytope*
\begin{pf}
    By \autoref{thm:gale-duality-bijection}, the $\Phi_{[n]\setminus [1]}$ designs of $J(n,k)$ are in bijection with the faces of the $\Phi_{[1]}$-eigenpolytope of $J(n, k)$.
    This eigenpolytope is described by Godsil \cite[Section 5]{godsil1998eigenpolytope}; its faces correspond to the elements of $\binom{[n]}{k}$ which contain $S$ and do not intersect $T$ for some choice of disjoint subsets $S, T$. The facets come from the case of $S = \emptyset, T = \set{i}$, or $S = \set{i}, T = \emptyset$. The minimal $\Phi_{[n]\setminus [1]}$-designs are the complements of these facets, and by \autoref{thm:erdos-ko-rado-designs} they are uniformly weighted. In the case $S = \emptyset, T = \set{i}$, the complementary minimal design is the star $\cD_i$, and in the case $S = \set{i}, T = \emptyset$, the complementary minimal design is the complement of $\cD_i$. Any face of the eigenpolytope is the intersection of some facets, thus any other design is the union of the minimal designs. 
\end{pf}

\subsubsection{Proof of \autoref{lem:scheme-to-graphical-design-cond}}
\lemschemetographicaldesigncond*
\begin{pf}
    If the hypothesis holds, then $G_i$ has $k+1$ distinct eigenspaces, given by the column spaces of $U_\ell$. Further, $U_0, U_1, \ldots, U_k$ are indexed in increasing order of eigenvalues of the Laplacian of $G_i$.
    Thus the columns of the $U_\ell$ for $\ell \in [t]$ together span the first $t$ eigenspaces of $G_i$ i.e. the space $\Phi_{[t]}$. Therefore, we can conclude that $\cD$ is a $[t]$-design of the scheme, if and only if $U_\ell^\top \bbone_\cD = 0$ for $\ell \in [t]$, if and only if $\cD$ is a $\Phi_{[t]}$-design of $G_i$. 
\end{pf}

\subsection{Proofs from \autoref{sec:S_n-section}}

\subsubsection{Proof of \autoref{thm:t-wise-uniform-designs}} \label{sec:proof-t-uniform-designs}

\twiseuniformdesigns*
\begin{pf}
    Consider two $t$-tuples $I = (i_1, \dots, i_t)$ and $J=(j_1, \dots, j_t)$ of distinct elements drawn from $[n]$. 
    We define the function $f_{I, J}: S_n \to \ZZ$ such that $f_{I, J}(\sigma) = 1$ if $\sigma(i_k) = j_k$ for all $k = 1, \dots, t$, and $f_{I, J}(\sigma) = 0$ otherwise. For any subset $\cD \subset S_n$, the count $\frac{1}{|\cD|}\sum_{\sigma \in \cD} f_{I, J}(\sigma)$ is the proportion of $\sigma\in \cD$ sending $I$ to $J$. Therefore $\cD$ is a $t$-wise uniform set precisely when this proportion for $\cD$ is the same as the proportion for all of $S_n$, for all choices of $I, J$. Specifically, 
    
    \[\frac{1}{|T|}\sum_{\sigma \in T} f_{I, J}(\sigma) = \frac{1}{|S_n|}\sum_{\sigma \in S_n} f_{I, J}(\sigma). \]

    So $\cD$ is a $t$-wise uniform set precisely when $\cD$ averages all functions $f_{I, J}$. Let $W$ be the space of functions $S_n \to \RR$ spanned by $f_{I, J}$, and let $\Phi$ be the space spanned by the functions in $\Phi_p$ with $p_1 \ge n-t$. Kuperberg, Lovett and Peled \cite[Section 3.4.2]{kuperberg2017combstructures} show that $W$ is the same as $\Phi$ (which they denote $W'$). Therefore averaging all $f_{I, J}$ is the same as averaging all of $\Phi$ i.e. $\cD$ is a $t$-wise uniform set if and only if it is a $\Phi$-design.
\end{pf}

\subsubsection{Proof of \autoref{prop:transposition-order}}

\proptranspositionorder*
\begin{pf}
    Let $\lambda'_p$ be the adjacency matrix eigenvalue of the eigenvector coming from the partition $p$. We can explicitly compute this eigenvalue from the parts of $p$ using \cite[Equation 3]{konstantinova2024transpositioneigenvalues}: if $p = (p_1, \dots, p_k) \vdash n$, we have
    \[\lambda'_p = \frac12\sum_{i=1}^k p_i(p_i - 2i + 1).\]

    Suppose $n$ is odd with $n = 2m+1$. Take $p = (m+1, 1, \dots, 1)$ with $m+1$ parts, and $q = (m, m, 1)$ with $3$ parts. Then we can compute
    \[\lambda'_p = \frac12\paren{(m+1)m + \sum_{i=2}^{m+1} (2 - 2i)} = \frac{m(m+1)}{2} - \paren{\sum_{i=2}^{m+1} i-1} = 0; \]
    \[\lambda'_q = \frac12\paren{m(m-1) + m(m-3) - 4} = m^2 - 2m - 2.\]
    Since $n \ge 6$ we have $m \ge 3$ and thus $m^2 - 2m - 2 > 0$. Therefore $\lambda'_p < \lambda'_q$. The Laplacian flips the order of eigenvalues compared to the adjacency matrix, so we have $\lambda_p > \lambda_q$. We also have $p_1 > q_1$, and thus $p$ and $q$ are the desired partitions.

    If $n$ is even with $n = 2m+2$, take $p = (m+1, 1, \dots, 1)$ with $m+2$ parts and $q = (m, m, 2)$ with $3$ parts. A similar computation as above shows that $\lambda'_p = - m - 1$ and $\lambda'_q = m^2 - 2m - 3$. If $n \ge 6$ we have $m \ge 2$ and thus $\lambda'_p \le \lambda'_q$ (with equality only if $m = 2$). Therefore in this case we also have $\lambda_p \ge \lambda_q$ and $p_1 > q_1$. 
\end{pf}

\subsubsection{Proof of \autoref{thm:symm-subgrp-designs}}
\thmsymmsubgrpdesigns*
\begin{pf}
    Translation by an element of $S_n$ is an automorphism of the Cayley graphs on $S_n$. Therefore, if we show that $S_t$, the subgroup of permutations fixing $[n]\setminus [t]$, is a design averaging $\Phi_p$ with $p_1 < n-t$, then all of its cosets $\cD = \sigma S_t \sigma'$ will also be designs averaging the same eigenvectors.
    
    The characteristic vector $\bbone_{S_t}$ on $S_n$ is the induced character of the trivial representation of $S_t$. This trivial representation corresponds to the partition $(t) \vdash t$. By the branching rule \cite[Chapter 7]{fulton1997tableaux}, the induced representation only has components along the irreducible representations of $S_n$ whose corresponding Young diagrams can be obtained by adding boxes to the diagram for $(t)$. Since adding boxes cannot reduce the first part, all the partitions $p$ with $p_1 < n-t$ have representations without components in the induced representation i.e. they have zero inner product with the induced character $\bbone_{S_t}$.Therefore $S_t$ averages the corresponding eigenvectors, and so it averages $\Phi_p$ with $p < n-t$.
\end{pf}

\subsubsection{Proof of \autoref{prop:birkhoff-eigenpolytope}}
\propbirkhoffeigenpolytope*
\begin{pf}
    The irreducible representation corresponding to the partition $(n-1, 1)$ is called the standard representation of $S_n$; let $\rho: S_n \to \RR^{(n-1)\times (n-1)}$ be the orthogonal standard representation. It is well known that $\rho(\sigma)$ is the matrix acting on $\RR^{n-1}$, by permuting the vertices of the $n$-simplex (of dimension $n-1$) according to the permutation $\sigma$. 
    The set $\Phi_{(n-1, 1)}$ contains $(n-1)^2$ functions $S_n \to \RR$, coming from the entries of the $\rho(\sigma)$ as $\sigma$ ranges over $S_n$. Taking the $(n-1)^2$ values for each $\sigma \in S_n$, we get a vector in $\RR^{(n-1)^2}$. The $\Phi_{(n-1, 1)}$-eigenpolytope is the convex hull of these vectors. 

    Now, the $n$-simplex also naturally lives in $\RR^n$ as the convex hull of the standard basis vectors $e_1, \dots, e_n$. This embedding is called the probability simplex. It lies in the affine hyperplane $\sum x_i = 1$, which is orthogonal to the all-ones vector $\bbone$.
   
    For any matrix $M \in \RR^{(n-1)\times(n-1)}$, we can define its lift $\tilde{M} \in \RR^{n\times n}$ as the matrix acting like $M$ when restricted to the affine hyperplane $\sum x_i = 1$ and satisfying $\tilde{M}\bbone = \bbone$. Since $\tilde{M}$ preserves all affine hyperplanes perpendicular to $\bbone$, we have that $\bbone^\top \tilde{M} v = \bbone^\top  v$ for all $v \in \RR^n$. Therefore $\bbone^\top \tilde{M} = \bbone^\top$ i.e. $\tilde{M}$ is doubly stochastic. The space of doubly stochastic matrices is isomorphic to $\RR^{(n-1)\times(n-1)}$, and thus the lifting operation is an affine isomorphism (i.e. a bijective affine tranformation with affine inverse). 

    The vertices of the probability simplex are the standard basis vectors. Thus the lifts $\widetilde{\rho(\sigma)}$ are the matrices in $\RR^{n\times n}$ permuting the standard basis vectors, which are precisely the permutation matrices. Therefore, their convex hull is the Birkhoff Polytope. Since hyperplanes and their orientations are preserved by affine isomorphisms, the original $\Phi_{(n-1, 1)}$-eigenpolytope before the lift must also be combinatorially isomorphic to the Birkhoff polytope.

    It is well known that the Birkhoff polytope has $n^2$ facets, one for each pair $i, j \in [n]$, comprised of the vertices corresponding to permutations with $\pi(i) \ne j$ \cite[Example 0.12]{ziegler1995polytopes}. By \autoref{thm:gale-duality-bijection}, the minimal designs averaging all eigenvectors outside $\Phi_{(n-1, 1)}$ are the complements of these subsets, containing the permutations with $\pi(i) = j$. These are precisely the cosets of $S_{n-1}$, and by \autoref{thm:symm-subgrp-designs} they are uniformly weighted. Since these are all the minimal designs, any other designs must be come from unions of such sets.
\end{pf}

\subsection{Proofs from \autoref{sec:mycielskian-section}}

\subsubsection{Proof of \autoref{thm:myc-central-vtx}}

\thmmyccentralvtx*
\begin{pf}
    We will argue that the eigenvectors averaged by the central vertex are those with eigenvalues of the form $\phi\mu_i$ and $\overline{\phi}\mu_i$.
    Suppose $i \in \{2, \dots, n\}$, and $x_i$ is an eigenvector of $A(G)$ with eigenvalue $\mu_i$. We claim that $[x_i\ -\overline{\phi} x_i\ 0]^\top$ is an eigenvector for $A(\cM(G))$ with eigenvalue $\phi\mu_i$. 
    Multiplying the adjacency matrix $A(\cM(G))$ with the vector $[x_i\ -\overline{\phi} x_i\ 0]^\top$, we get
    \[
        \begin{bmatrix}
            A & A & 0\\
            A & 0 & \bbone \\
            0 & \bbone^\top  & 0
        \end{bmatrix}\begin{bmatrix}
            x_i \\ -\overline\phi x_i \\ 0
        \end{bmatrix} = \begin{bmatrix}
            A(x_i - \overline\phi x_i) \\ Ax_i \\ -\overline\phi \bbone^\top x_i
        \end{bmatrix}.
    \]
    Since $\phi, \overline\phi$ are the roots of $x^2 - x - 1$, we have $1 - \overline\phi = \phi$. Therefore $A(x_i - \overline\phi x_i) = A(\phi x_i) = (\phi\mu_i) x_i$. Also we have $1 = -\phi\overline\phi$, and thus $Ax_i = \mu_ix_i = (\phi\mu_i)(-\overline\phi x_i)$. Finally, since $G$ is $d$-regular, it has $\bbone$ as an adjacency eigenvector with eigenvalue $d$. Thus we have $\bbone^\top x_i = 0$. From the above computation, we can conclude that 
    \[
        \begin{bmatrix}
            A(x_i - \overline\phi x_i) \\ Ax_i \\ -\overline\phi \bbone^\top x_i
        \end{bmatrix} = \begin{bmatrix} (\phi\mu_i) x_i \\ (\phi\mu_i)(-\overline\phi x_i) \\ 0 \end{bmatrix}.
    \]
    Therefore $[x_i\ -\overline{\phi} x_i\ 0]^\top$ is an eigenvector of $A(\cM(G))$ with eigenvalue $\phi\mu_i$. A similar computation shows that $[x_i\ -\phi x_i\ 0]^\top$ is an eigenvector for $A(\cM(G))$ with eigenvalue $\overline\phi\mu_i$. Note that all $2n-2$ such eigenvectors have value 0 at $u$. 
        
    If an adjacency eigenvector is orthogonal to $\bbone$, a similar argument to the Laplacian case (\autoref{sec:laplacian-preliminaries}) shows that a design averages the eigenvector if its values at the design vertices sum to zero. Using the fact that $\bbone^\top x_i = 0$, a routine check shows that both the eigenvectors $[x_i\ -\overline{\phi} x_i\ 0]^\top$ and $[x_i\ -\phi x_i\ 0]^\top$ are orthogonal to the all-ones vector in $\RR^{2n+1}$. 
    Therefore all $2n-2$ eigenvectors averaging the eigenvalues $\phi\mu_i$ and $\overline\phi\mu_i$ for $i = 2, \dots, n$ are orthogonal to $\bbone$  and have value zero at the central vertex $u$. We conclude that the design $\set{u}$ averages all these eigenvectors. 
\end{pf}

\subsubsection{Proof of \autoref{prop:myc-designs-from-original-graph}}

\propmycdesignsfromoriginalgraph*
\begin{pf}
    Note that $V$ and $V'$ are both copies of the vertex set of the original graph, and thus we can identify the sets $\cD \cap V$ and $\cD \cap V'$ as designs of $G$. Suppose $x_i$ is an eigenvector of $A(G)$ with eigenvalue $\mu_i$. For $v' \in V'$, define $x_i(v') = x_i(v)$ where $v$ is the copy of $v'$ in the original vertex set $V$. Then to show that $\cD\cap V$ and $\cD \cap V'$ are designs of $G$ averaging $x_i$, we must show 
    \[\sum_{v\in\cD\cap V} x_i(v) = \sum_{v' \in \cD \cap V'} x_i(v') = 0.\]
    From the proof of \autoref{thm:myc-central-vtx}, if $x_i$ is an eigenvector of $A(G)$ with eigenvalue $\mu_i$, the eigenvectors of $A(\cM(G))$ with eigenvalue $\phi\mu_i$ and $\overline\phi\mu_i$ are $[x_i\ -\overline{\phi} x_i\ 0]^\top$ and $[x_i\ -\phi x_i\ 0]^\top$ respectively. 
    If $\cD$ is a design averaging both, then the values of the eigenvectors at the vertices of $\cD$ sum to zero, which implies
    \[\sum_{v \in \cD \cap V} x_i(v) - \overline\phi\sum_{v' \in \cD \cap {V'}} x_i(v') = 0 = \sum_{v \in \cD \cap V} x_i(v) - \phi\sum_{v' \in \cD \cap {V'}} x_i(v').\]
    Solving the linear system, we have $\sum_{v \in \cD \cap V} x_i(v) = \sum_{v' \in \cD \cap {V'}} x_i(v') = 0$. Thus both $\cD \cap V$ and $\cD \cap {V'}$ are designs of $G$ averaging $x_i$.
\end{pf}

\subsubsection{Proof of \autoref{prop:myc-design-sizes}}

We first establish a lemma.
\begin{lem}[label=lem:rational-weights-lemma]
    Consider a graph $G = (V, E)$, with an adjacency eigenvalue $\alpha \notin \QQ$. Then $\alpha$ will be an algebraic integer with minimal polynomial $P(x)$ over $\QQ$. Any design $\cD \subset V$ which averages the eigenspace for $\alpha$ will also average the eigenspaces for its conjugates $\alpha'$ (the other roots of $P(x)$).
\end{lem}
\begin{pf}
    Let $K$ be the splitting field of $\alpha$; this is the smallest field in which all the roots of the minimal polynomial $P(x)$ are contained. For any conjugate $\alpha'$ of $\alpha$, we have a field automorphism $\sigma$ of $K / \QQ$ with $\sigma(\alpha) = \alpha'$. 

    If $\phi_\alpha$ is an eigenvector of $A$ with eigenvalue $\alpha$, we have $(A - \alpha I)\phi_\alpha = 0.$ We apply the automorphism $\sigma$ to this equation; note that $A, I$ have rational entries and are thus fixed by $\sigma$. Then we have
    \[(\sigma(A) - \sigma(\alpha)\sigma(I))\sigma(\phi_\alpha) = (A - \alpha'I)\sigma(\phi_\alpha) = 0.\]
    Therefore, $\sigma(\phi_\alpha)$ is an eigenvector of $A$ with eigenvalue $\alpha'$ i.e. $\sigma$ is a bijection between the eigenspaces of $\alpha$ and $\alpha'$. Now, since $\cD$ averages the eigenspace for $\alpha$,
    \[ \sum_{v \in \cD} \phi_\alpha(v) =  \frac{1}{|V|} \sum_{v \in V} \phi_\alpha(v).\]        
    Applying $\sigma$ to the equation, we get
    \[ \sum_{v \in \cD} \sigma(\phi_\alpha(v)) =  \frac{1}{|V|} \sum_{v \in V} \sigma(\phi_\alpha(v)).\]        
    Therefore $\cD$ also averages $\sigma(\phi_\alpha)$, and thus averages the eigenspace for $\alpha'$.
\end{pf}

In the setting of weighted designs (\autoref{sec:weighted-designs-and-eigenpolytopes}), \autoref{lem:rational-weights-lemma} generalizes to hold for all designs with rational weights.

\propmycdesignsizes*
\begin{pf}
    Any design $\cD$ averaging the first $2n-2$ eigenvectors in Mycielskian order must average all the eigenvectors with eigenvalues $\phi\mu_i$ and $\overline\phi\mu_i$ for $i = 2, \dots, n$. By \autoref{prop:myc-designs-from-original-graph}, $\cD \cap V$ and $\cD \cap {V'}$ are designs of the original graph averaging the eigenvectors with eigenvalue $\mu_i$ for $i = 2, \dots, n$. Now, the eigenvector with eigenvalue $\mu_1 = d$ is constant, and thus it is averaged by all designs. Therefore $\cD \cap V$ and $\cD \cap {V'}$ are designs averaging all eigenvectors of $A(G)$; such a design must contain all vertices or be empty. 
    If $\cD$ is not $\set{u}$, either $\cD \cap V$ or $\cD \cap {V'}$ is nonempty and then the respective intersection is all of $V$ or $V'$. We conclude that either $\cD = \set{u}$ or it contains at least $n$ vertices.

    The remaining three eigenvectors of the Mycielskian have eigenvalues coming from the roots of $t^3 - dt^2 - (n+d^2)t + dn$ i.e. from conjugate eigenvalues. By \autoref{lem:rational-weights-lemma}, any design averaging one of these eigenvectors must average all three. This shows that designs averaging more than $2n-2$ eigenvectors in Mycielskian order must average all of them, and thus those designs must contain all vertices.
\end{pf}

\subsection*{Acknowledgements} 
The authors gratefully acknowledge helpful conversations with William Monty McGovern, Thomas Rothvoss, Isaiah Siegl and Varun Shah.

\bigskip
\printbibliography

@book{godsil-meagher, place={Cambridge}, series={Cambridge Studies in Advanced Mathematics}, title={Erdõs–Ko–Rado Theorems: Algebraic Approaches}, publisher={Cambridge University Press}, author={Godsil, Christopher and Meagher, Karen}, year={2015}, collection={Cambridge Studies in Advanced Mathematics}}

@book{godsil-algcomb,
author = {Godsil, C.},
year = {1993},
title = {Algebraic Combinatorics (1st ed.).},
publisher = {Routledge},
doi = {https://doi.org/10.1201/9781315137131}
}

@article{johnsen1985lineare,
  title={Lineare Abh{\"a}ngigkeiten von Einheitswurzeln.},
  author={Johnsen, Karsten and Albrechts, Ch},
  journal={Elemente der Mathematik},
  volume={40},
  pages={57--59},
  year={1985}
}

@article{pillai1933arithmetic,
  title={On an arithmetic function},
  author={Pillai, Subbayya Sivasankaranarayana},
  journal={Annamalai University Journal},
  volume={2},
  pages={243--248},
  year={1933}
}

@book {alon2008probabilistic,
    AUTHOR = {Alon, Noga and Spencer, Joel H.},
     TITLE = {The probabilistic method},
   EDITION = {Third},
      NOTE = {With an appendix on the life and work of Paul Erd\H os},
 PUBLISHER = {John Wiley \& Sons, Inc., Hoboken, NJ},
      YEAR = {2008},
     PAGES = {xviii+352},
      ISBN = {978-0-470-17020-5},
   MRCLASS = {60-02 (05C80 60C05 60F99 60G42)},
  MRNUMBER = {2437651},
       DOI = {10.1002/9780470277331},
       URL = {https://doi.org/10.1002/9780470277331},
}

@book {macwilliams-sloane,
    AUTHOR = {MacWilliams, F. J. and Sloane, N. J. A.},
     TITLE = {The theory of error-correcting codes.},
    SERIES = {North-Holland Mathematical Library},
    VOLUME = {16},
 PUBLISHER = {North-Holland Publishing Co.},
      YEAR = {1977},
      ISBN = {0-444-85009-0},
   MRCLASS = {94A10},
  MRNUMBER = {465509},
MRREVIEWER = {Ian\ Blake},
}

@article{babecki2021codes,
  title={Codes, cubes, and graphical designs},
  author={Babecki, Catherine},
  journal={Journal of Fourier Analysis and Applications},
  volume={27},
  number={5},
  pages={81},
  year={2021},
  publisher={Springer},
  url={https://link.springer.com/article/10.1007/s00041-021-09852-z}
}

@article{babecki2023galeduality,
  title={Graphical designs and gale duality},
  author={Babecki, Catherine and Thomas, Rekha R},
  journal={Mathematical Programming},
  volume={200},
  number={2},
  pages={703--737},
  year={2023},
  publisher={Springer},
  url={https://link.springer.com/article/10.1007/s10107-022-01861-0}
}

@article {babecki2022whatis,
    AUTHOR = {Babecki, Catherine},
     TITLE = {What is{$\ldots$} a graphical design?},
   JOURNAL = {Notices Amer. Math. Soc.},
  FJOURNAL = {Notices of the American Mathematical Society},
    VOLUME = {69},
      YEAR = {2022},
    NUMBER = {9},
     PAGES = {1571--1573},
      ISSN = {0002-9920,1088-9477},
   MRCLASS = {05E18 (05C10 05C25)},
  MRNUMBER = {4501025},
}

@article{steinerberger2020designs,
  title={Generalized designs on graphs: sampling, spectra, symmetries},
  author={Steinerberger, Stefan},
  journal={Journal of graph theory},
  volume={93},
  number={2},
  pages={253--267},
  year={2020},
  publisher={Wiley Online Library},
  url={https://arxiv.org/pdf/1803.02235}
}

@article {babecki2024universality,
    AUTHOR = {Babecki, Catherine and Shiroma, David},
     TITLE = {Eigenpolytope universality and graphical designs},
   JOURNAL = {SIAM J. Discrete Math.},
  FJOURNAL = {SIAM Journal on Discrete Mathematics},
    VOLUME = {38},
      YEAR = {2024},
    NUMBER = {1},
     PAGES = {947--964},
      ISSN = {0895-4801,1095-7146},
   MRCLASS = {05C50 (52B12 68Q17 68R10 90C57)},
  MRNUMBER = {4712817},
MRREVIEWER = {Hong\ Zhang},
       DOI = {10.1137/22M1528768},
       URL = {https://doi.org/10.1137/22M1528768},
}

@article {babecki2024sparsedesigns,
    AUTHOR = {Al-Thani, Hessa and Babecki, Catherine and Mart\'inez Mori, J.
              Carlos},
     TITLE = {Sparse graphical designs via linear programming},
   JOURNAL = {Oper. Res. Lett.},
  FJOURNAL = {Operations Research Letters},
    VOLUME = {56},
      YEAR = {2024},
     PAGES = {Paper No. 107145, 7},
      ISSN = {0167-6377,1872-7468},
   MRCLASS = {05C50 (90C05 90C35)},
  MRNUMBER = {4781608},
       DOI = {10.1016/j.orl.2024.107145},
       URL = {https://doi.org/10.1016/j.orl.2024.107145},
}

@article {zhu2023bch,
    AUTHOR = {Zhu, Yan},
     TITLE = {Optimal and extremal graphical designs on regular graphs
              associated with classical parameters},
   JOURNAL = {Des. Codes Cryptogr.},
  FJOURNAL = {Designs, Codes and Cryptography. An International Journal},
    VOLUME = {91},
      YEAR = {2023},
    NUMBER = {8},
     PAGES = {2737--2754},
      ISSN = {0925-1022,1573-7586},
   MRCLASS = {05B30 (05E30 94B15)},
  MRNUMBER = {4618186},
       DOI = {10.1007/s10623-023-01231-7},
       URL = {https://doi.org/10.1007/s10623-023-01231-7},
}

@article {kuperberg2017combstructures,
    AUTHOR = {Kuperberg, Greg and Lovett, Shachar and Peled, Ron},
     TITLE = {Probabilistic existence of regular combinatorial structures},
   JOURNAL = {Geom. Funct. Anal.},
  FJOURNAL = {Geometric and Functional Analysis},
    VOLUME = {27},
      YEAR = {2017},
    NUMBER = {4},
     PAGES = {919--972},
      ISSN = {1016-443X,1420-8970},
   MRCLASS = {05D40 (05A15 05B15 05B30 05C65)},
  MRNUMBER = {3678505},
       DOI = {10.1007/s00039-017-0416-9},
       URL = {https://doi.org/10.1007/s00039-017-0416-9},
}

@article {konstantinova2024transpositioneigenvalues,
    AUTHOR = {Konstantinova, Elena V. and Kravchuk, Artem},
     TITLE = {Distinct eigenvalues of the transposition graph},
   JOURNAL = {Linear Algebra Appl.},
  FJOURNAL = {Linear Algebra and its Applications},
    VOLUME = {690},
      YEAR = {2024},
     PAGES = {132--141},
      ISSN = {0024-3795,1873-1856},
   MRCLASS = {05C25 (05E10)},
  MRNUMBER = {4722337},
       DOI = {10.1016/j.laa.2024.03.011},
       URL = {https://doi.org/10.1016/j.laa.2024.03.011},
}

@article {balakrishnan2012mycielskian,
    AUTHOR = {Balakrishnan, R. and Kavaskar, T. and So, Wasin},
     TITLE = {The energy of the {M}ycielskian of a regular graph},
   JOURNAL = {Australas. J. Combin.},
  FJOURNAL = {The Australasian Journal of Combinatorics},
    VOLUME = {52},
      YEAR = {2012},
     PAGES = {163--171},
      ISSN = {1034-4942,2202-3518},
   MRCLASS = {05C50 (05C90)},
  MRNUMBER = {2917924},
}

@article{mycielski1955,
  author = {Mycielski, J.},
  journal = {Colloq. Math},
  number = {3},
  title = {Sur le coloriage des graphes},
  volume = {161-162},
  year = {1955}
}

@book {chung1997spectral,
    AUTHOR = {Chung, Fan R. K.},
     TITLE = {Spectral graph theory},
    SERIES = {CBMS Regional Conference Series in Mathematics},
    VOLUME = {92},
 PUBLISHER = {Conference Board of the Mathematical Sciences, Washington, DC;
              by the American Mathematical Society, Providence, RI},
      YEAR = {1997},
     PAGES = {xii+207},
      ISBN = {0-8218-0315-8},
   MRCLASS = {58G99 (05C50 35P05 46N20 47N20)},
  MRNUMBER = {1421568},
MRREVIEWER = {Robert\ Brooks},
}

@article {goethals1977sphericaldesigns,
    AUTHOR = {Delsarte, P. and Goethals, J. M. and Seidel, J. J.},
     TITLE = {Spherical codes and designs},
   JOURNAL = {Geometriae Dedicata},
  FJOURNAL = {Geometriae Dedicata},
    VOLUME = {6},
      YEAR = {1977},
    NUMBER = {3},
     PAGES = {363--388},
   MRCLASS = {05B99},
  MRNUMBER = {485471},
MRREVIEWER = {Michel\ Deza},
       DOI = {10.1007/bf03187604},
       URL = {https://doi.org/10.1007/bf03187604},
}

@article {godsil1998eigenpolytope,
    AUTHOR = {Godsil, C. D.},
     TITLE = {Eigenpolytopes of distance regular graphs},
   JOURNAL = {Canad. J. Math.},
  FJOURNAL = {Canadian Journal of Mathematics. Journal Canadien de
              Math\'ematiques},
    VOLUME = {50},
      YEAR = {1998},
    NUMBER = {4},
     PAGES = {739--755},
      ISSN = {0008-414X,1496-4279},
   MRCLASS = {05E30 (05C50)},
  MRNUMBER = {1638611},
MRREVIEWER = {Jack\ H.\ Koolen},
       DOI = {10.4153/CJM-1998-040-8},
       URL = {https://doi.org/10.4153/CJM-1998-040-8},
}

@article {deza1977frankl,
    AUTHOR = {Frankl, P\'eter and Deza, Mikhail},
     TITLE = {On the maximum number of permutations with given maximal or
              minimal distance},
   JOURNAL = {J. Combinatorial Theory Ser. A},
  FJOURNAL = {Journal of Combinatorial Theory. Series A},
    VOLUME = {22},
      YEAR = {1977},
    NUMBER = {3},
     PAGES = {352--360},
      ISSN = {0097-3165},
   MRCLASS = {05A10},
  MRNUMBER = {439648},
MRREVIEWER = {Luc\ Teirlinck},
       URL =
              {http://www.sciencedirect.com/science/article/pii/0097316577900097},
}

@article {erdos1961intersection,
    AUTHOR = {Erd\H os, P. and Ko, Chao and Rado, R.},
     TITLE = {Intersection theorems for systems of finite sets},
   JOURNAL = {Quart. J. Math. Oxford Ser. (2)},
  FJOURNAL = {The Quarterly Journal of Mathematics. Oxford. Second Series},
    VOLUME = {12},
      YEAR = {1961},
     PAGES = {313--320},
      ISSN = {0033-5606,1464-3847},
   MRCLASS = {04.60},
  MRNUMBER = {140419},
MRREVIEWER = {S.\ Ginsburg},
       DOI = {10.1093/qmath/12.1.313},
       URL = {https://doi.org/10.1093/qmath/12.1.313},
}

@article {golubev2020extremal,
    AUTHOR = {Golubev, Konstantin},
     TITLE = {Graphical designs and extremal combinatorics},
   JOURNAL = {Linear Algebra Appl.},
  FJOURNAL = {Linear Algebra and its Applications},
    VOLUME = {604},
      YEAR = {2020},
     PAGES = {490--506},
      ISSN = {0024-3795,1873-1856},
   MRCLASS = {05B99 (05C35 05C50 05C69 05C70 35J05 35P05 35R02)},
  MRNUMBER = {4123763},
       DOI = {10.1016/j.laa.2020.07.012},
       URL = {https://doi.org/10.1016/j.laa.2020.07.012},
}

@book {grigoryan2018intro,
    AUTHOR = {Grigor'yan, Alexander},
     TITLE = {Introduction to analysis on graphs},
    SERIES = {University Lecture Series},
    VOLUME = {71},
 PUBLISHER = {American Mathematical Society, Providence, RI},
      YEAR = {2018},
     PAGES = {viii+150},
      ISBN = {978-1-4704-4397-9},
   MRCLASS = {60J10 (05C25 05C50 05C63 05C81 58J35)},
  MRNUMBER = {3822363},
MRREVIEWER = {Serguei\ Popov},
       DOI = {10.1090/ulect/071},
       URL = {https://doi.org/10.1090/ulect/071},
}

@article {linderman2020numerint,
    AUTHOR = {Linderman, George C. and Steinerberger, Stefan},
     TITLE = {Numerical integration on graphs: where to sample and how to
              weigh},
   JOURNAL = {Math. Comp.},
  FJOURNAL = {Mathematics of Computation},
    VOLUME = {89},
      YEAR = {2020},
    NUMBER = {324},
     PAGES = {1933--1952},
      ISSN = {0025-5718,1088-6842},
   MRCLASS = {05C50 (05C70 35J05 35P05 41A55 65D30)},
  MRNUMBER = {4081923},
MRREVIEWER = {Farideh\ Heydari},
       DOI = {10.1090/mcom/3515},
       URL = {https://doi.org/10.1090/mcom/3515},
}

@book {nica2018intro,
    AUTHOR = {Nica, Bogdan},
     TITLE = {A brief introduction to spectral graph theory},
    SERIES = {EMS Textbooks in Mathematics},
 PUBLISHER = {European Mathematical Society (EMS), Z\"urich},
      YEAR = {2018},
     PAGES = {viii+156},
      ISBN = {978-3-03719-188-0},
   MRCLASS = {05-01 (05C50 11T24 15-01 20C15)},
  MRNUMBER = {3821579},
       DOI = {10.4171/188},
       URL = {https://doi.org/10.4171/188},
}

@article {sobolev1962cubature,
    AUTHOR = {Sobolev, S. L.},
     TITLE = {Cubature formulas on the sphere which are invariant under
              transformations of finite rotation groups},
   JOURNAL = {Dokl. Akad. Nauk SSSR},
  FJOURNAL = {Doklady Akademii Nauk SSSR},
    VOLUME = {146},
      YEAR = {1962},
     PAGES = {310--313},
      ISSN = {0002-3264},
   MRCLASS = {65.55},
  MRNUMBER = {141225},
MRREVIEWER = {A.\ H.\ Stroud},
}

@article {steinerberger2021spherical,
    AUTHOR = {Steinerberger, Stefan},
     TITLE = {Spectral limitations of quadrature rules and generalized
              spherical designs},
   JOURNAL = {Int. Math. Res. Not. IMRN},
  FJOURNAL = {International Mathematics Research Notices. IMRN},
      YEAR = {2021},
    NUMBER = {16},
     PAGES = {12265--12280},
      ISSN = {1073-7928,1687-0247},
   MRCLASS = {58J50 (31C15 41A55)},
  MRNUMBER = {4300225},
MRREVIEWER = {Patrick\ Guidotti},
       DOI = {10.1093/imrn/rnz176},
       URL = {https://doi.org/10.1093/imrn/rnz176},
}

@article {steinerberger2025equidistribution,
    AUTHOR = {Steinerberger, Stefan and Thomas, Rekha R.},
     TITLE = {Random walks, equidistribution and graphical designs},
   JOURNAL = {Adv. in Appl. Math.},
  FJOURNAL = {Advances in Applied Mathematics},
    VOLUME = {165},
      YEAR = {2025},
     PAGES = {Paper No. 102837, 11},
      ISSN = {0196-8858,1090-2074},
   MRCLASS = {05C48 (05C81)},
  MRNUMBER = {4848395},
       DOI = {10.1016/j.aam.2024.102837},
       URL = {https://doi.org/10.1016/j.aam.2024.102837},
}

@article {kharaghani2004hadamard,
    AUTHOR = {Kharaghani, H. and Tayfeh-Rezaie, B.},
     TITLE = {A {H}adamard matrix of order 428},
   JOURNAL = {J. Combin. Des.},
  FJOURNAL = {Journal of Combinatorial Designs},
    VOLUME = {13},
      YEAR = {2005},
    NUMBER = {6},
     PAGES = {435--440},
      ISSN = {1063-8539,1520-6610},
   MRCLASS = {05B20},
  MRNUMBER = {2221851},
MRREVIEWER = {Arne\ Winterhof},
       DOI = {10.1002/jcd.20043},
       URL = {https://doi.org/10.1002/jcd.20043},
}

@article {dokovic2014hadamard,
    AUTHOR = {\DJ okovi\'c, Dragomir \v Z. and Golubitsky, Oleg and
              Kotsireas, Ilias S.},
     TITLE = {Some new orders of {H}adamard and skew-{H}adamard matrices},
   JOURNAL = {J. Combin. Des.},
  FJOURNAL = {Journal of Combinatorial Designs},
    VOLUME = {22},
      YEAR = {2014},
    NUMBER = {6},
     PAGES = {270--277},
      ISSN = {1063-8539,1520-6610},
   MRCLASS = {05B20 (05B10)},
  MRNUMBER = {3193950},
MRREVIEWER = {Mieko\ Yamada},
       DOI = {10.1002/jcd.21358},
       URL = {https://doi.org/10.1002/jcd.21358},
}

@book {ziegler1995polytopes,
    AUTHOR = {Ziegler, G\"unter M.},
     TITLE = {Lectures on polytopes},
    SERIES = {Graduate Texts in Mathematics},
    VOLUME = {152},
 PUBLISHER = {Springer-Verlag, New York},
      YEAR = {1995},
     PAGES = {x+370},
      ISBN = {0-387-94365-X},
   MRCLASS = {52Bxx},
  MRNUMBER = {1311028},
MRREVIEWER = {Margaret\ M.\ Bayer},
       DOI = {10.1007/978-1-4613-8431-1},
       URL = {https://doi.org/10.1007/978-1-4613-8431-1},
}

@book {fulton1997tableaux,
  title={Young Tableaux: With Applications to Representation Theory and Geometry},
  author={Fulton, W.},
  isbn={9780521567244},
  lccn={95047484},
  series={London Mathematical Society Student Texts},
  year={1997},
  publisher={Cambridge University Press}
}

@article {ellis2011intersecting,
    AUTHOR = {Ellis, David and Friedgut, Ehud and Pilpel, Haran},
     TITLE = {Intersecting families of permutations},
   JOURNAL = {J. Amer. Math. Soc.},
  FJOURNAL = {Journal of the American Mathematical Society},
    VOLUME = {24},
      YEAR = {2011},
    NUMBER = {3},
     PAGES = {649--682},
      ISSN = {0894-0347,1088-6834},
   MRCLASS = {05D10 (05A05 05C50 05E10 20C30)},
  MRNUMBER = {2784326},
MRREVIEWER = {Norihide\ Tokushige},
       DOI = {10.1090/S0894-0347-2011-00690-5},
       URL = {https://doi.org/10.1090/S0894-0347-2011-00690-5},
}

@article {filmus2017comment,
      title={A comment on Intersecting Families of Permutations}, 
      author={Yuval Filmus},
      year={2017},
      eprint={1706.10146},
      archivePrefix={arXiv},
      primaryClass={math.CO},
      url={https://arxiv.org/abs/1706.10146}, 
}

@article {wilson1984erdoskorado,
    AUTHOR = {Wilson, Richard M.},
     TITLE = {The exact bound in the {E}rd\H os-{K}o-{R}ado theorem},
   JOURNAL = {Combinatorica},
  FJOURNAL = {Combinatorica. An International Journal of the J\'anos Bolyai
              Mathematical Society},
    VOLUME = {4},
      YEAR = {1984},
    NUMBER = {2-3},
     PAGES = {247--257},
      ISSN = {0209-9683},
   MRCLASS = {05A05 (05C35 05C65)},
  MRNUMBER = {771733},
MRREVIEWER = {Noga\ Alon},
       DOI = {10.1007/BF02579226},
       URL = {https://doi.org/10.1007/BF02579226},
}

@article {delsarte1973thesis,
    AUTHOR = {Delsarte, P.},
     TITLE = {An algebraic approach to the association schemes of coding
              theory},
   JOURNAL = {Philips Res. Rep. Suppl.},
  FJOURNAL = {Philips Research Reports. Supplements},
      YEAR = {1973},
    NUMBER = {10},
     PAGES = {vi+97},
   MRCLASS = {94A10},
  MRNUMBER = {384310},
MRREVIEWER = {Ian\ Blake},
}

\end{document}